\documentclass[11pt]{amsart}
\usepackage[margin=1in]{geometry}
\usepackage{amsmath,amssymb,mathtools}
\usepackage{mathrsfs}
\mathtoolsset{showonlyrefs}
\usepackage{mathtools}
\usepackage{microtype}
\usepackage[svgnames]{xcolor}
\usepackage[colorlinks,linkcolor=FireBrick,citecolor=RoyalBlue,urlcolor=RoyalBlue]{hyperref}

\usepackage{tikz}
\usetikzlibrary{arrows.meta}

\numberwithin{equation}{section}
\newtheorem{theorem}{Theorem}[section]

\newtheorem{lemma}[theorem]{Lemma}

\theoremstyle{definition}
\newtheorem{definition}[theorem]{Definition}
\theoremstyle{remark}
\newtheorem{remark}[theorem]{Remark}
\newcommand{\E}{\mathbb E}
\newcommand{\Pp}{\mathbb P}
\newcommand{\one}{\mathbf1}
\newcommand{\tr}{\operatorname{tr}}
\newcommand{\cH}{\mathcal H}
\newcommand{\cT}{\mathcal T}
\newcommand{\cF}{\mathcal F}
\newcommand{\cP}{\mathcal P}
\newcommand{\cR}{\mathcal R}

\title[]{A second eigenvalue lower bound for the directed configuration model}
\author{Charles Bordenave}
\address{I2M, CNRS, Aix-Marseille Universit\'e}
\email{charles.bordenave@univ-amu.fr}
\author{Yizhe Zhu}
\address{Department of Mathematics, University of Southern California}
\email{yizhezhu@usc.edu}
\date{}

\begin{document}
\begin{abstract}
We prove a high-probability lower bound for the second eigenvalue in modulus of the transition matrix of the simple random walk on a directed configuration model with minimal out-degree at least two and maximal degree not too large. The lower bound matches the upper bound found by Coste (2021) on the same model. Our proof bypasses the need to establish convergence of the empirical eigenvalue distribution. Instead, it reduces the spectral lower bound to a central limit theorem for the trace of the transition matrix at powers much larger than the logarithmic scale, a range that previous results failed to reach.  Both the general strategy and the central limit theorem might be of independent interest.
\end{abstract}
\maketitle

\section{Introduction}

\subsection*{Motivation from random matrix theory}

Consider a sequence of possibly random matrices $A \in M_n(\mathbb C)$ of growing dimension $n \to \infty$. In many situations, we are interested in the asymptotics as $n \to \infty$ of the {\em spectral radius} of $A$ defined as
$$
\varrho(A) = \max \{ |\lambda| : \lambda \in \sigma(A)\}.
$$
The spectral radius is upper bounded by the operator norm: 
$$
\varrho(A) \leq \| A \|_{\rm op} = \sup_{ x \ne 0} \frac{ \| A x\|_2 }{\| x \|_2}.
$$
In particular, for any integer $k \geq 1$, we have $\varrho(A)^k = \varrho(A^k) \leq \| A^k \|_{\rm op}$. Gelfand's theorem asserts that this upper bound is asymptotically sharp in the sense that 
$$
\varrho(A) = \lim_{k \to \infty} \| A^k \|_{\rm op}^{1/k}.
$$
This limit is, however, not uniform over the dimension $n$, and this formula may not be very useful in high dimensions.

The {\em trace method} is a classical general strategy to upper bound the spectral radius of random matrices; see the monographs \cite{AndersonGuionnetZeitouni2010,BaiSilverstein2010,Tao2012}. A foundational application to the extreme eigenvalues of
symmetric random matrices is due to F\"uredi and Koml\'os
\cite{FurediKomlos1981}. It starts with the following chain of inequalities: 
$$
\varrho(A)^{2k} =  \varrho(A^k)^2 \leq \| A ^k  \|_{\rm op}^{2}  = \| A^k (A^k)^* \|_{\rm op} \leq \tr \left( A^k  (A^k)^*  \right).
$$ 
Then, we may expand the trace as a sum of products of matrix entries and, for many sequences of random matrices $A$, we may derive a probabilistic upper bound on the growth rate of trace powers of the form: $\tr \left( A^k  (A^k)^*  \right) \leq n (1+\epsilon_k)^{2k} \rho_+^{2k}$ for some $\rho_+ >0$ and where $\epsilon_k$ goes to $0$ as $k \to \infty$. Therefore, we find that $\varrho (A) \leq n^{1/2k} (1+\epsilon_k)\rho_+ $. If this holds for $k \gg \log n$,  we then deduce  that $\varrho (A) \leq  (1+O(\max( \epsilon_k, \log n / k))) \rho_+$.  

For example, if $A = X/ \sqrt n$ with $(X_{ij})$ independent, centered, variance $1$ and uniformly bounded support, it follows from Geman \cite{zbMATH03978046}, one of the earliest appearances of this method, that for some $c_1 >0$, $\mathbb E \tr (A^k (A^k)^* ) \leq n k^{c_1}$ for $k \leq n^{1/c_1}$, leading to $\varrho(A) \leq  1 + O(n^{-c_2})$ with high probability for some $c_2 >0$. Since then, there have been numerous improvements and further examples of such successful computations in random matrix theory; see, for example \cite{Soshnikov1999,Vu2007,BandeiraVanHandel2016,BenaychGeorgesBordenaveKnowles2020}.

If  $A$ is a normal matrix (that is if $A A^* = A^* \hspace{-2pt}A$), then the same strategy gives also a lower bound. Indeed, for normal matrices, we have $\rho(A)= \| A \|_{\rm op}$ and thus 
$$
\varrho(A)^{2k}  = \| (AA^*)^ {k}  \|_{\rm op} \geq \frac 1 n \tr \left( (AA^*)^{k}  \right).
$$
In particular, any lower bound of the form $\tr \left( (AA^*)^ {k}  \right) \geq n (1-\epsilon_k )^{2k}\rho_-^{2k}$ now gives $\varrho(A) \geq (1-\epsilon_k) \rho_-$. Note that for the lower bound, there is no need to let $k \to \infty$ jointly with $n \to \infty$. If $\rho_+ = \rho_-$, then we have obtained an asymptotic equivalent to $\varrho(A)$ with a controlled error.

In the non-normal case, however, it is not straightforward to extend this strategy to get a lower bound. The usual route takes the following detour. We first prove that the empirical distribution of the eigenvalues of $A$ converges weakly to a compactly supported measure $\mu$ on $\mathbb C$. Then, we observe that it implies that $\liminf_{n \to \infty} \varrho(A) \geq \rho_- $ with $\rho_- = \sup \{ |\lambda| : \lambda \in \mathrm{supp} (\mu) \}$.

For example, if $A = X/ \sqrt n$ with $(X_{ij})$ iid centered entries of variance $1$, then, in probability, the empirical distribution of eigenvalues converges weakly to the uniform measure on the unit disc; see \cite{zbMATH06415652} for a survey. Together with the upper bound $\varrho(A)\leq 1+o(1)$ from \cite{bordenave2022convergence}, we thus get in probability, $\lim_n \varrho (A) = 1$ in this case.

Proving the weak convergence of the empirical distribution of eigenvalues is usually not an easy task, and one might wonder if a more direct and quantitative approach exists. One goal of this work is to illustrate one such possible approach. It still relies on the trace powers, but now it exploits its {\em trace fluctuations for large powers}.  More precisely, we write: 
$$
\varrho(A) \geq \frac{1}{n}\sum_{i=1}^n |\lambda_i(A)| \geq \left|  \frac{1}{n}\sum_{i=1}^n \lambda_i(A) \right| = \left| \frac 1 n \tr (A) \right|. 
$$
Hence for any integer $k \geq 1$, 
\begin{equation}\label{eq:radiusLB00}
\varrho(A) \geq   \left| \frac 1 n  \tr (A^k) \right|^{1/k}. 
\end{equation}
For fixed $k$, at first order in $n$, this lower bound may be crude. For example, if $\varrho (A) = O(1)$ and if the empirical distribution of eigenvalues of $A$ converges to a radial measure $\mu$ on $\mathbb C$ then for any integer $k \geq 1$,
$$
\lim_{n \to \infty} \frac{1}n \tr (A^ k) =  \int_{\mathbb C} \lambda^k d\mu(\lambda) = 0.
$$
As we discussed above, this case notably happens if $A = X/ \sqrt n$ with $(X_{ij})$ i.i.d. centered entries of variance $1$. On the other hand, for this random matrix, under a proper moment condition on $X_{ij}$, for any fixed $k \geq 1$, it is well-known \cite{RiderSilverstein2006,NourdinPeccati2010} that, as $n \to \infty$,
$$
\frac{\tr (A^k)}{\sqrt k} \Longrightarrow Z ,
$$
where $Z$ is a Gaussian variable of unit variance. Moreover, central limit theorems for linear eigenvalue
statistics have been established for analytic test functions
in real elliptic ensembles \cite{ORourkeRenfrew2016}, and for
sufficiently smooth test functions in non-Hermitian matrices with i.i.d. entries
\cite{CipolloniErdosSchroder2023,CipolloniErdosSchroder2021}.

More quantitatively, for any random matrix model $A$, assume that we can prove that weakly, along a sequence $k = k_n$ with $k \gg \log n$,
$$
\frac{\tr (A^k)}{(1+\epsilon_k)^k \rho_-^k} \Longrightarrow Z,
$$
where $\rho_- >0$ and $Z$ is a random variable such that $\mathbb P ( Z = 0) = 0 $. Then, $|Z|^{1/k} \geq 1 - O(\log n/k)$ with high probability. We thus get from \eqref{eq:radiusLB00} that, with high probability,
$$
 \varrho(A) \geq (1 + \epsilon_k) \left( 1 - O \left( \frac{\log n }{k} \right)\right)  \rho_-. 
$$
Therefore, the fluctuations of trace powers provide a lower bound on the spectral radius.

To our knowledge, this strategy has not been used so far. In the present work, we follow this strategy for the study of the spectral gap of the simple random walk on a random directed graph sampled according to the configuration model, which we define below.

\subsection*{Model and motivation from spectral theory of Markov chains}
For integer $n \geq 1$, prescribe non-negative integers $(d_i^+,d_i^-)_{i=1}^n$ such that
\begin{equation}\label{eq:degrees}
 \sum_i d_i^+=\sum_i d_i^-=N.
\end{equation}

If $N \geq 1$, we consider a random directed multigraph $G$ on $[n]$ sampled from the configuration model. It is defined as follows. Attach $d_i^+$
outgoing half-edges, called \emph{heads}, and $d_i^-$ incoming
half-edges, called \emph{tails}, to vertex $i$. Write $E^+(i)$
and $E^-(i)$ for these two sets, and use bold letters
$\mathbf e$ and $\mathbf f$ for heads and tails, respectively. The unions
\[
 E^+=\bigsqcup_{i=1}^n E^+(i),\qquad
 E^-=\bigsqcup_{i=1}^n E^-(i)
\]
both have cardinality $N$. 
Choose arbitrary labels $[N]$ for each
of $E^+$ and $E^-$, sample a uniform permutation $\sigma\in\mathcal S_N$,
and let $\sigma:E^+\to E^-$ also denote the induced bijection: the
head with label $a$ is matched to the tail with label $\sigma(a)$. In other words, $\sigma:E^+\to E^-$ is a uniform bijection.
If $\mathbf e\in E^+(i)$ and $\sigma(\mathbf e)\in E^-(j)$,
we place one directed edge $i\to j$.

Write $G$ for the resulting directed configuration multigraph;
loops and parallel edges are allowed. Assuming $\min_i d_i^+ \geq 1$, for $i,j\in[n]$, the $(i,j)$-entry of the transition probability matrix of $G$ is
\begin{equation}
\label{eq:P^k}
 P(i,j)=\frac{\#\{\mathbf e\in E^+(i):
                  \sigma(\mathbf e)\in E^-(j)\}}{d_i^+}.
\end{equation}

We  order the eigenvalues of $P$
by modulus, including multiplicity: 
\[
 \lambda_1=1,\qquad
 |\lambda_2|\ge\cdots\ge|\lambda_n|.
\]

The matrix $P$ is the Markov transition kernel of the simple random walk on $G$. It is non-reversible on a directed graph. The spectral gap of $P$ is $1 - |\lambda_2|$.  It measures the asymptotic mixing rate to equilibrium; we refer to  \cite{Norris1997,MR3726904}. For example, if $P$ is aperiodic and irreducible, it admits a unique invariant probability measure $\pi$ and then, for any probability measure $\mu$ on $[n]$, 
$$
\lim_{t \to \infty} \| \mu P^t - \pi \|_{\mathrm{TV}} ^{1/t} \le |\lambda_2|,
$$
where for a signed measure $\nu$, its total variation is $\| \nu \|_{\mathrm{TV}} = (1/2) \cdot \sum_{i} | \nu(i)| $. 

The invariant measure $\pi$ of $P$ defined in \eqref{eq:P^k} and its mixing time were studied in \cite{bordenave2018random}. The study of the spectrum of $P$ was initiated by Coste~\cite{coste2021spectral}, where an upper bound on $|\lambda_2|$ was computed using the trace method outlined above. In this work, we prove a matching lower bound for $|\lambda_2|$; see Remark \ref{rk:simon} below.

\subsection*{Main results} Define
\begin{equation}\label{eq:rho}
\Delta = \max_{i} \{ d_i^+ , d_i^-\}, \qquad \delta=\min_i d_i^+ \quad \hbox{ and } \quad \rho=\sqrt{\frac1 N\sum_{i=1}^n \frac{d_i^-}{d_i^+}}.
\end{equation}
We assume that for all sufficiently large integers $n$,
\begin{equation}\label{eq:gap}
\delta \geq 2 \quad \hbox{ and } \quad \delta\rho\ge 1. 
 \end{equation}

 The condition $\delta\rho\ge 1$ is satisfied when the out-degree sequence is concentrated.
For example, since $\rho \geq 1/ \sqrt \Delta$, the condition  $\delta\rho\ge 1$ is implied by the assumption $\delta^2 \geq \Delta$.  Similarly, let $\kappa = \min_i d_i^- / d_i^+ \in [0,1]$ and $d = N/n$ be the average in/out-degree. Then $\rho \geq \sqrt{\kappa / d}$ and the condition $\delta\rho\ge 1$ is implied by the assumption $\kappa \delta^2 \geq d $.

We finally consider another parameter of the degree sequence. Set 
$$
\beta_3 = \frac1 N\sum_{i=1}^n \frac{d_i^-}{(d_i^+)^2}
$$
and consider some $\theta \geq 0$ such that 
\begin{equation}\label{eq:defbeta3}
\beta_3 \leq e^{-\theta} \rho ^3. 
\end{equation}
Indeed, if \eqref{eq:gap} holds, then $1/d_i^+ \leq 1/\delta \leq \rho$ and thus \eqref{eq:defbeta3} holds for $\theta = 0$. More is true: we will show  in Lemma \ref{lem:beta} that $\theta  = (\delta - \sqrt{\delta})/\Delta>0 $ always satisfies \eqref{eq:defbeta3}. If $G$ is out-degree regular, that is if $d_i^+ = d$ for all $i$, then $\rho = 1/ \sqrt d$, $\beta_3 = 1/d^2$ and thus we may take $\theta =  (\log d )/2$.

We prove the following lower bound on the second eigenvalue of $P$ provided that $\theta$ in \eqref{eq:defbeta3} is not too small and that $\Delta$ is not too large.

\begin{theorem}\label{thm:main}
Assume  \eqref{eq:gap}-\eqref{eq:defbeta3} with $\log \Delta    =  o( \frac{\log n}{\log\log n})$ and $\theta \geq \sqrt{\log \Delta} /(\log n )^{1+o(1)}$. 
We have
\begin{equation}\label{eq:lower}
 \Pp\{|\lambda_2|\ge(1-\varepsilon_n)\rho\}\longrightarrow 1,
\end{equation}
for any sequence $\varepsilon_n \to 0$ such that \begin{equation}\label{eq:epsilon}
 \varepsilon_n \gg  \log \Delta  \cdot\frac{  \log\log n}{ \log n }. 
\end{equation}
\end{theorem}
As pointed out below \eqref{eq:defbeta3}, since $\theta \geq  ( 2 - \sqrt{2} )/ \Delta$ in the worst case, the condition on $\theta$ in Theorem~\ref{thm:main} is automatically satisfied  if $\Delta \leq (\log n)^{1+o(1)}$ and \eqref{eq:gap} holds.

\begin{remark}[A matching upper bound]\label{rk:simon}
    When $\Delta= O(1)$, Coste~\cite{coste2021spectral} proves that,
with high probability,
$|\lambda_2|\leq\max\{\rho,1/\underline \delta\}+o(1)$, where
$\underline \delta=\min\{d_i^+,d_i^- :  i\in [n]\}$.
A close inspection of  the proof in \cite{coste2021spectral} shows that 
$\underline \delta$ can be replaced by the minimum out-degree $\delta$.
The same argument therefore yields the sharper bound
$|\lambda_2|\leq\max\{\rho,1/\delta\}+o(1)$.
Under assumption~\eqref{eq:gap}, this upper bound reduces to
$\rho+o(1)$ and matches our lower bound.
\end{remark}

As explained above, we prove Theorem~\ref{thm:main} by establishing the following central limit theorem for the traces of powers of $P$. This is our main technical contribution.

\begin{theorem}\label{lem:CLT}Assume  \eqref{eq:gap}-\eqref{eq:defbeta3} with $\Delta = n^{o(1)}$ and some $\theta >0$. Let $t = t_n$ be an integer sequence such that 
\begin{equation}\label{eq:parametert}
1 \ll t  \ll \frac{ (\log n )^2}{ \log (\Delta)  \log\log n} \quad \hbox{ and } \quad  \theta^{-2+o(1)} \leq t.
\end{equation}
Then, we have
    \begin{align}\label{eq:CLT_P}
     \frac{\tr P^{t}-1}{\sqrt{t}\,\rho^{t}}
 \Longrightarrow\mathcal N(0,1).
\end{align}
\end{theorem}

\begin{proof}[Proof of Theorem~\ref{thm:main} with Theorem~\ref{lem:CLT}] 
Since the conclusion becomes stronger when $\varepsilon_n$
is decreased, we may assume without loss of generality that
\begin{align}\label{eq:epsilon_upper_bound}
 \varepsilon_n\le
 \frac{\log\Delta\,(\log\log n)^2}{\log n}.
\end{align}
The trace identity gives 
\[
 |\tr P^{t}-1|
 =\left|\sum_{i=2}^n\lambda_i^{t}\right|
 \le(n-1)|\lambda_2|^{t}.
\]
Take $t =   \lfloor 2 ( \log n  ) / \varepsilon_n \rfloor $ where $\varepsilon_n$ satisfies \eqref{eq:epsilon}, \eqref{eq:epsilon_upper_bound}, and $\varepsilon_n=o(1)$. Then the two conditions in  \eqref{eq:parametert} hold. Also, we have, for sufficiently large $n$,
\[
 \varepsilon_n t \geq \log n.
\]
Hence, for sufficiently large $n$,
$(1-\varepsilon_n)^{t}
 \le e^{-\varepsilon_n t}\le n^{-1}$.
Therefore
\begin{align}
    \mathbb P (|\lambda_2|\leq (1-\varepsilon_n) \rho) &\leq\mathbb P \left( \left|\frac{\tr P^{t} -1}{\sqrt{t} \rho^{t}}\right| \leq \frac{n-1}{\sqrt{t}}(1-\varepsilon_n)^{t}\right)\\
    &\leq \mathbb P \left( \left|\frac{\tr P^{t} -1}{\sqrt{t} \rho^{t}}\right| \leq \frac{1}{\sqrt{t}}\right)\leq \mathbb P \left( \left|\frac{\tr P^{t} -1}{\sqrt{t} \rho^{t}}\right| \leq \epsilon\right)
\end{align}
for any $\epsilon>0$ when $n$ is sufficiently large.
For every fixed $\epsilon>0$, Theorem \ref{lem:CLT}
therefore gives
\[
 \limsup_{n\to\infty}
 \Pp\{|\lambda_2|\le(1-\varepsilon_n)\rho\}
 \le\Pp\{|Z|\le\epsilon\},
 \qquad Z\sim\mathcal N(0,1).
\]
Letting $\epsilon\downarrow0$, and using that the Gaussian
distribution has no atom at zero, we prove \eqref{eq:lower}.
\end{proof}

 Central limit theorems for traces and cycle counts at lengths
up to the diameter scale  $\log_{d-1}(n)$ have been established for random
$d$-regular graphs and digraphs 
in~\cite{zbMATH06206998,johnson2011exchangeable,coste2024characteristic},
using Poisson approximations of cycle counts and Gaussian limits for Poisson
random variables with diverging means. Beyond the graph diameter scale, Garmo \cite{garmo1999asymptotic} obtained weighted Poisson limits for normalized cycle count fluctuations in fixed-degree random regular graphs.
Theorem~\ref {lem:CLT} holds with $t\gg\log_{\Delta} n$, beyond the range covered by previous
results. As we can see in the proof above, reaching this scale is essential for
obtaining the sharp spectral lower bound via~\eqref{eq:radiusLB00}.

\subsection*{Proof overview} A natural approach to proving Theorem~\ref{lem:CLT} is by the moment method. However, as in the proofs of Alon's second eigenvalue conjecture \cite{friedman2008proof,bordenave2015new}, tangles (see Definition~\ref{def:tangle}) can dominate the moments when $t\gg \log n$. To illustrate this obstruction, consider the directed configuration
model with $d_i^{\pm}=d$ for a fixed $d\ge2$, so that
$\rho=d^{-1/2}$. Let $\Omega$
be the event that vertices $1,2$ are each isolated with $d$
self-loops. Counting matchings gives for some constant $C_d$ depending only on $d$,
\[
 \Pp(\Omega)
 =(d!)^2\frac{(nd-2d)!}{(nd)!}
 \ge C_d n^{-2d}.
\] On $\Omega$, we have $\tr P^t\ge2$ and  when $t\gg \log n$,
\begin{align}
    \E\left( \frac{\tr P^t-1}{\sqrt{t}\rho^t}\right)^2\geq \frac{d^t}{t}\mathbb P(\Omega) \geq C_d\frac{d^{t}}{tn^{2d}}\to\infty.
\end{align}
Thus, when $t\gg \log n$ in Theorem~\ref{lem:CLT}, the normalized trace converges in
distribution to a standard Gaussian, even though its second moment
diverges. The isolated vertices with $d$ self-loops are examples
of the tangles creating this issue. 

To overcome this obstruction, we use the tangle-free path
decomposition of \cite{bordenave2015new,coste2021spectral}
and apply the moment method to a modified trace. Pairs of traces following the same directed cycle
of length $t$ give the Gaussian moments, while all other paths have relatively negligible total contribution. On a high-probability tangle-free event, the modified trace agrees with
$\tr P^t-1$ up to a negligible scalar correction after normalization. This transfers the Gaussian
limit to the original trace and completes the proof.

\subsection*{Obstruction to lower bounds on $|\lambda_2|$}

As stressed above, the Markov chain $P$ is non-reversible, which in linear algebra terms means that $P \ne P^*$ where the adjoint is taken with respect to the invariant measure. A consequence of non-reversibility is that there are no known simple lower bounds on $|\lambda_2|$; see \cite{MehtaSchulman2020} for relations between graph expansion and spectral gap.

The failure of a deterministic Alon-Boppana bound for $d$-regular digraphs is illustrated by the de Bruijn digraphs.  For integer $s\geq1$, the digraph $DB(d,s)$ has vertex set 
$
 V(DB(d,s))=[d]^s
$
and edge set 
$
 E(DB(d,s))=\left\{
 (a_1,\ldots,a_s)\longrightarrow(a_2,\ldots,a_s,t):
 a_i,t\in[d]
 \right\}.
$
It is an aperiodic $d$-regular digraph on $d^s$ vertices, but all of its
nontrivial eigenvalues vanish
\cite{delorme1998spectrum,parzanchevski2020ramanujan}.

\subsection*{Further relations to previous works and extensions}

As mentioned above, the simple random walk on the directed configuration model was previously studied in  \cite{bordenave2018random,coste2021spectral}. Some of the techniques developed in this work build upon a refinement of the trace method developed in \cite{bordenave2015new,coste2021spectral}.

Despite the lack of general Alon-Boppana-type lower bounds for $|\lambda_2|$, there have been some recent works on such lower bounds for adjacency or non-backtracking operators of sparse random digraphs based on limiting spectral distributions \cite{SahSahasrabudheSawhney2023,he2026oriented,mckenzie2026alon,coste2026spectral,basak2018circular,litvak2021circular}.

Notably, a special case of the directed configuration model is the random $d$-regular digraph (when $d_{i}^{+}=d_{i}^{-}=d$ for all $i\in [n]$). It was recently proved in \cite{he2026oriented} that, in probability, the empirical distribution of eigenvalues of $P$ converges weakly to the oriented Kesten-McKay law. This implies $|\lambda_2(P)|\geq  1/\sqrt{d}-o(1)$ with high probability for random $d$-regular digraphs. Our present work gives a quantitative estimate on this lower bound.

There are several other closely related sparse random matrices which could also be studied with our techniques. The simpler model of this form is the anisotropic random walk on a directed random Schreier graph:  that is, $A = \sum_{i=1}^d p_i S_i$ where $S_i$ are independent uniform permutation matrices of size $n$ and $\sum_i p_i = 1$, $p_i \geq 0$. It follows from \cite{BordenaveCollins2019,BordenaveQiuZhang2020} that $|\lambda_2(A)| \leq (1+o(1)) \|p\|_2$. An easy adaptation of the present work gives the matching lower bound.

More generally, in \cite{BordenaveQiuZhang2020}
and  \cite{Timhadjelt2023}, the authors study random sparse bistochastic matrices  of the form $QS$ and $p Q + (1-p)S$, respectively,  where $S$ is a uniform permutation matrix, $Q$ is a bistochastic matrix, and $0< p <1$. They establish sharp upper bounds on $|\lambda_2|$. It would be interesting to investigate the corresponding lower bounds with the techniques presented in this work.

\section{Proof of the main result}
It remains to prove Theorem~\ref{lem:CLT}. Our proof uses the
tangle-free trace method developed in
\cite{bordenave2015new,coste2021spectral}.

\subsection{Estimates on a random matching}

For a random permutation from $\mathcal S_N$, define a matching matrix $M$ such that
\[
 M_{\mathbf e\mathbf f}=\one_{\{\sigma(\mathbf e)=\mathbf f\}},
 \qquad \underline M_{\mathbf e\mathbf f}
 =M_{\mathbf e\mathbf f}-N^{-1}.
\]
Consider a sequence of half-edges
$\gamma=(\gamma_1,\ldots,\gamma_{2k})$, with
$\gamma_{2t-1}\in E^+$ and $\gamma_{2t}\in E^-$.
 The distinct oriented edges
in the sequence form the set
\[
 E_\gamma=\bigl\{\{\gamma_{2t-1},\gamma_{2t}\}:1\le t\le k\bigr\}.
\]

Since $E^+$ and $E^-$ are disjoint, each element of $E_\gamma$
specifies one head and one tail. We call such a pair an \emph{edge};
it need not belong to the sampled matching.
For $e=\{\mathbf e,\mathbf f\}\in E_\gamma$, its
\emph{multiplicity} is
\[
 \#\{t\in[k]:\gamma_{2t-1}=\mathbf e,\ \gamma_{2t}=\mathbf f\}.
\]
The edge is called \emph{single} when this number is one. It is
\emph{consistent} when
\[
 e\cap e'=\varnothing \quad  \text{ for any } \quad e'\in E_\gamma,\ e'\ne e.
\]
Distinct edges sharing a head or a tail cannot both belong to the same
matching, and we call them \textit{inconsistent}.

A factor $M_{\mathbf e\mathbf f}$
is called \emph{raw}, and a factor $\underline M_{\mathbf e\mathbf f}$
is called \emph{centered}. We have the following estimate for the expectation of a mixed product for edges in a random matching. 

\begin{lemma}\label{lem:matching}
Let $1\le k\le\sqrt N$ and $0\le k_0\le k$. Let
$\gamma=(\gamma_1,\ldots,\gamma_{2k})$ be a sequence of half-edges
such that $\gamma_{2t-1}\in E^+$ and $\gamma_{2t}\in E^-$
for every $1\le t\le k$.
Put $a=|E_\gamma|$. Among the edges occurring at positions
$t\le k_0$, let $a_1$ count the consistent single edges and
$b$ the inconsistent single edges, with multiplicities taken in
the whole sequence $\gamma$. When $k_0=0$, put $a_1=b=0$ and
interpret the empty product as one. We have for some constant $C>0$,
\begin{equation}\label{eq:entry}
 \left|\E
 \prod_{t=1}^{k_0}\underline M_{\gamma_{2t-1}\gamma_{2t}}
 \prod_{t=k_0+1}^{k}M_{\gamma_{2t-1}\gamma_{2t}}\right|
 \le \frac{C\,2^b}{ N^{a}}\left(\frac{3k}{\sqrt N}\right)^{a_1}.
\end{equation}
When $k_0=k$, assume all edges are consistent, and every edge  has multiplicity
 at least two, then
\begin{equation}\label{eq:entry-leading}
\left|  \E\prod_{t=1}^{k}\underline M_{\gamma_{2t-1}\gamma_{2t}} -\frac{1}{(N)_a}\right|\leq \frac{2k}{N\cdot (N)_a}.
\end{equation}
where  $(N)_a=N(N-1)\cdots(N-a+1)$.
\end{lemma}
\begin{proof}
For $k_0\ge1$, \eqref{eq:entry} is
\cite[Proposition 28]{bordenave2015new}, see also \cite{ALCO_2026__9_2_465_0} for the connection to Weingarten calculus. For $k_0=0$, an
inconsistent product is zero, whereas a consistent product has
expectation $(N)_a^{-1}$. Since $a\le k\le\sqrt N$,
\begin{align}\label{eq:Na_ratio}
 0\le\log\frac{N^a}{(N)_a}
 =\sum_{h=0}^{a-1}-\log(1-h/N) \leq \sum_{h=0}^{a-1}\frac{h}{N-a+1}
 \le\frac{a(a-1)}{2(N-a+1)}\le\frac12,
\end{align}
so \eqref{eq:entry} holds in this case as well.
For \eqref{eq:entry-leading}, let $h_1,\ldots,h_a\ge2$ be the
multiplicities, so $\sum_{i=1}^a h_i=k$.
For a matching indicator $M_{\mathbf e\mathbf f}$,
\[
 (M_{\mathbf e\mathbf f}-N^{-1})^h
 =\alpha_h M_{\mathbf e\mathbf f}+(-N^{-1})^h, \quad \text{where}
 \quad 
 \alpha_h=(1-N^{-1})^h-(-N^{-1})^h.
\]
For $N\ge2$, we have $\alpha_h\in[0,1]$ and
$1-\alpha_h\le h/N$. The following elementary inequality holds:
\begin{equation}\label{eq:elementary}
    0\le1-\prod_{i=1}^a\alpha_{h_i}
 \le\sum_{i=1}^a(1-\alpha_{h_i})\le\frac{k}{N}. 
\end{equation}
Let $\{\mathbf{e}_t,\mathbf{f}_t\}, 1\leq t\leq a$ be the distinct edges. Then
\begin{align}
    \E\left[\prod_{t=1}^{a}\underline M_{\mathbf{e}_t,\mathbf{f}_t}^{h_t}\right] -\frac{1}{(N)_a}=\E \prod_{t=1}^{a}(\alpha_{h_t} M_{\mathbf {e}_t\mathbf {f}_t}+(-N^{-1})^{h_t})-\frac{1}{(N)_a}.
\end{align}
In the product $\E \prod_{t=1}^{a}(\alpha_{h_t} M_{\mathbf {e}_t\mathbf {f}_t}+(-N^{-1})^{h_t})$, the term selecting all $a$ indicator random variables contributes
$\frac{1}{(N)_a}\prod_{i=1}^a\alpha_{h_i}$, and from \eqref{eq:elementary},
\[
\left| \frac{1}{(N)_a}\prod_{i=1}^a\alpha_{h_i}
 -\frac{1}{(N)_a}\right|\le\frac{k}{N (N)_a}.
\]
Since $h_t\geq 2$, each term taking $(a-j)$ indicators and $j$ scalars contributes at most
\[\frac{1}{(N)_{a-j}}\cdot N^{-2j}\leq \frac{1}{(N)_a} N^{-j}\] in absolute  value. The total approximation error is at most 
\begin{align}
    \left|  \E\prod_{t=1}^{k}\underline M_{\gamma_{2t-1}\gamma_{2t}} -\frac{1}{(N)_a}\right|&\leq \frac{k}{N (N)_a}+\frac{1}{(N)_a} \sum_{j=1}^a \binom{a}{j} N^{-j}\\
    &=\frac{1}{(N)_a}\left( \frac{k}{N} +(1+N^{-1})^a-1\right)\leq \frac{2k}{N\cdot  (N)_a},
\end{align}
where in the last inequality, we use $a\leq k/2$.
\end{proof}

\subsection{Tangle-freeness and telescoping decomposition}
The matching estimate Lemma \ref{lem:matching} applies to arbitrary half-edge sequences. To
use it for powers of $P$, we now impose the vertex adjacency condition
that makes such a sequence a path.

\begin{definition}[Path]\label{def:path}
For $t\ge1$, a path of length $t$ from $u$ to $v$ is a sequence
$\gamma=(\mathbf e_1,\mathbf f_1,\ldots,\mathbf e_t,\mathbf f_t)$
for which there are vertices $v_0,\ldots,v_t$ satisfying
\[
 v_0=u,\qquad v_t=v,\qquad
 \mathbf e_s\in E^+(v_{s-1}),\quad
 \mathbf f_s\in E^-(v_s),\quad 1\le s\le t.
\]
Vertices and half-edges may repeat. The path is defined before
sampling $\sigma$; it is \emph{realized} when
$\sigma(\mathbf e_s)=\mathbf f_s$ for every $s$.
Write $\cP^t(u,v)$ for the set of all these potential paths,
without requiring their edges to be realized or consistent.
\end{definition}

The \emph{support} $G(\gamma)$ is the directed multigraph with
vertex set $\{v_0,\ldots,v_t\}$ and one edge $v_{s-1}\to v_s$
for each distinct pair $(\mathbf e_s,\mathbf f_s)$. Repeated
uses of the same pair give only one support edge; distinct pairs
with the same endpoints give parallel edges. Loops are allowed.
Connectivity, cycles, and degrees
are computed after forgetting orientations but keeping all loops
and parallel edges. In particular, a loop contributes two to
the degree, and one to each of the in- and out-degrees.

\begin{definition}[Genus and tangle-free paths]\label{def:tangle}
For a connected multigraph $H$, its \emph{genus} is
\[
 g(H)=|E(H)|-|V(H)|+1.
\]
The graph $H$ is
\emph{tangle-free} if $g(H)\le1$ (a tree or a connected graph
with one cycle), and is a \emph{tangle} if $g(H)\ge2$.
A path $\gamma$ is tangle-free when its support $G(\gamma)$
is tangle-free, and \emph{tangled} otherwise. 
For integer $t \geq 1$, write $\cT^t(u,v)\subseteq\cP^t(u,v)$ for the tangle-free paths
of length $t$ from $u$ to $v$.
\end{definition}

Let $d_{\mathbf e}^+$ denote the out-degree of the vertex carrying
the head $\mathbf e$. Define the degree weight of a path
$\gamma=(\mathbf e_s,\mathbf f_s)_{s=1}^t$ by
\[
 \omega(\gamma)=\prod_{s=1}^t(d_{\mathbf e_s}^+)^{-1}.
\]
 The ordinary and
centered weights of a single step are
\[
 A(\mathbf e,\mathbf f)=\frac{M_{\mathbf e\mathbf f}}{d_{\mathbf e}^+},
 \qquad
 \underline A(\mathbf e,\mathbf f)
 =\frac{\underline M_{\mathbf e\mathbf f}}{d_{\mathbf e}^+}.
\]
The usual path expansion gives
\[
 P^t(u,v)=\sum_{\gamma\in\cP^t(u,v)}
              \prod_{s=1}^t A(\mathbf e_s,\mathbf f_s).
\]
Restricting the sum to tangle-free paths defines two matrices $P^{(t)}$ and $\underline P^{(t)}$ such that 
\begin{equation}\label{eq:tf}
 \begin{aligned}
 P^{(t)}(u,v)
 &=\sum_{\gamma\in\cT^t(u,v)}
       \prod_{s=1}^t A(\mathbf e_s,\mathbf f_s),\\
 \underline P^{(t)}(u,v)
 &=\sum_{\gamma\in\cT^t(u,v)}
       \prod_{s=1}^t\underline A(\mathbf e_s,\mathbf f_s).
 \end{aligned}
\end{equation}
Set $P^{(0)}=\underline{P}^{(0)}=I_n$.

Let $B^+(u,\ell)$ be the subgraph of $G$ induced by the vertices reachable from $u$
by a directed path of length at most $\ell$. The graph $G$ is
\textit{$\ell$-tangle-free} if every $B^+(u,\ell)$ is tangle-free.
Let $\cH_\ell$ be the event that $G$ is $\ell$-tangle-free. Then 
\[ P^{\ell}=P^{(\ell)} \qquad \text{on} \qquad \cH_\ell.\]

The following lemma extends the local tangle-free estimate in
\cite[Proposition 2.3]{coste2021spectral} to growing $\Delta$,
using the exploration argument of
\cite[Section 3.2]{bordenave2018random}.
\begin{lemma}[Tangle-freeness]\label{lem:local}
If $\Delta=n^{o(1)}$ and
$\ell \leq \lceil\alpha\log_\Delta n\rceil$ for a fixed $0<\alpha<1/4$,
then $G$ is $\ell$-tangle-free with high probability.
\end{lemma}
\begin{proof}
Starting with a fixed vertex, explore the
reachable vertices in order of directed distance, revealing all
outgoing heads at vertices of distance at most $\ell$, including
those at distance $\ell$.
This reveals every edge of its $\ell$-neighborhood and uses at most
\[
 \Delta\sum_{j=0}^\ell\Delta^j\le2\Delta^{\ell+1}
\]
matching pairs. Since $2\Delta^{\ell+1}\leq n^{\alpha+o(1)}=o(N)$, this number is
at most $N/2$ for sufficiently large $n$. A collision occurs when the chosen tail is
attached to an already discovered vertex. A tangle
requires at least two collisions.
The collision estimate in \cite[Lemma 5]{bordenave2018random},
followed by a union bound over the $n$ initial vertices, gives
\[
 \begin{aligned}
 \Pp\{G\text{ is not }\ell\text{-tangle-free}\}
 &\le n\,\frac{2\Delta^2(2\Delta^{\ell+1})^4}{(N/2)^2}\le\frac{32\Delta^{4\ell+6}}{n}
 \leq n^{4\alpha-1+o(1)}.
 \end{aligned}
\]
Since $0 < \alpha < 1/4$, the lemma follows.
\end{proof}

For $1 \leq k \leq \ell$, let $\cR^{\ell,k}(u,v)$
be the set of tangled paths whose first $(k-1)$ steps and
the last $(\ell-k)$ steps are separately tangle-free.
The $k$-th edge joins these two tangle-free paths. Define
\begin{equation}\label{eq:remainder}
 R^{\ell,k}(u,v)=
 \sum_{\gamma\in\cR^{\ell,k}(u,v)}
 \left(\prod_{s<k}A(\mathbf e_s,\mathbf f_s)\right)
 \frac1{d_{\mathbf e_k}^+}
 \left(\prod_{s>k}\underline A(\mathbf e_s,\mathbf f_s)\right).
\end{equation}

Write \[\pi^-(i)=\frac{d_i^{-}}{N}, \quad i\in [n]\] for the probability that a uniformly
chosen tail is attached to $i$. 
\begin{lemma}[Telescoping decomposition, Lemma 3.1 in \cite{coste2021spectral}]\label{lem:telescope}
For $\ell \geq 1$, the following identity holds:
\begin{equation}\label{eq:telescope}
 P^{(\ell)}=\underline P^{(\ell)}
 +\sum_{k=1}^\ell P^{(k-1)}\one(\pi^-)^\top\underline P^{(\ell-k)}
 -\frac1 N\sum_{k=1}^\ell R^{\ell,k}.
\end{equation}
\end{lemma}

Define, for integers $\ell , t \geq 1$,
\begin{equation}\label{eq:D}
 Q^{(\ell)}=\underline P^{(\ell)}
               -\frac1N\sum_{k=1}^{\ell}R^{\ell,k},
 \qquad z_\ell=\frac1n\one^\top Q^{(\ell)}\one,  \qquad z_{\ell,t}=z_{\ell}^{m-1}z_{\ell_0},
 \qquad X_{\ell, t}=\tr\bigl((Q^{(\ell)})^{m-1} Q^{(\ell_0)}\bigr),
\end{equation}
where $t = (m-1) \ell  + \ell_0$ and $m = \lceil t / \ell \rceil$.
The following lemma connects $\tr P^{t}$ with $(X_{\ell,t},z_{\ell,t})$ on the event $\mathcal H_{\ell}$.
\begin{lemma}\label{lem:identification}
 Let $\ell, t \geq 1$  be integers. On the tangle-free event $\cH_{\ell}$, 
    \begin{align}\label{eq:trace_identity}
    \tr P^{t}-1=   X_{\ell,t}-z_{\ell,t}.
    \end{align}
\end{lemma}
\begin{proof}
    On the event $\cH_{\ell}$, \eqref{eq:telescope} becomes
    \begin{align}\label{eq:free_telescope}
     P^{\ell}=\underline P^{(\ell)}
 +\sum_{k=1}^{\ell}P^{k-1}\one(\pi^-)^\top\underline P^{(\ell-k)}
 -\frac{1}{N}\sum_{k=1}^{\ell}R^{\ell,k}.
    \end{align}
Since $P$ is row stochastic,  $P^{k-1}\one=\one$ for any $k\geq 1$. From \eqref{eq:free_telescope}, we can write $P^{\ell}$ as a rank-one perturbation of $Q^{(\ell)}$:
\[
 P^\ell=Q^{(\ell)}+\one v_\ell^\top,
 \qquad \text{where} \qquad
 v_\ell^\top=\sum_{k=1}^{\ell}(\pi^-)^\top
                         \underline P^{(\ell-k)}.
\]
A simple calculation shows  $z_\ell=\frac1n\one^\top Q^{(\ell)}\one
       =1-v_\ell^\top\one$ and $
 Q^{(\ell)}\one=(1-v_\ell^\top\one)\one=z_{\ell}\one$.
If $m \geq 2$, we now telescope the matrix powers:
\[
 \begin{aligned}
 P^{(m-1)\ell}-(Q^{(\ell)})^{m-1}
 &=\sum_{j=1}^{m-1}(Q^{(\ell)})^{j-1}
            (P^\ell-Q^{(\ell)})P^{(m-1-j)\ell}\\
            &=\sum_{j=1}^{m-1}(Q^{(\ell)})^{j-1}
            \one v_\ell^\top P^{(m-1-j)\ell}=\sum_{j=1}^{m-1}z_\ell^{j-1}\one v_\ell^\top P^{(m-1-j)\ell}.
 \end{aligned}
\]
We next compose by $P^{\ell_0} = Q^{(\ell_0)} + \one v_{\ell_0}^{\top}$ (since $\ell_0 \le \ell$ implies $\cH_{\ell} \subset \cH_{\ell_0}$). We find
$$
P^t - (Q^{(\ell)})^{m-1} Q^{(\ell_0)} =\sum_{j=1}^{m-1}z_\ell^{j-1}\one v_\ell^\top P^{(m-1-j)\ell+\ell_0} + z_{\ell}^{m-1} \one v_{\ell_0}^\top.
$$
Taking traces and using $P^{k}\one=\one$  gives
\[
 \begin{aligned}
 \tr P^{t}-X_{\ell,t}
 &=\sum_{j=1}^{m-1}z_\ell^{j-1} v_\ell^\top \one + z_\ell^{m-1} v_{\ell_0}^\top \one =\sum_{j=1}^{m-1}z_\ell^{j-1}(1-z_\ell)
 + z_\ell^{m-1}  (1 - z_{\ell_0}) =1-z_\ell^{m-1} z_{\ell_0}.
 \end{aligned}
\]
This proves \eqref{eq:trace_identity}.
\end{proof}

\subsection{Trace power expansion}

For the remainder of this section, we set
\begin{equation}\label{eq:parameters}
 \ell=\left\lfloor\frac{\log n}{8\log\Delta}\right\rfloor,
\end{equation}
and consider an integer $t$ as in \eqref{eq:parametert}. We set $m = \lceil t / \ell \rceil$ and write $t$ as 
$$
t = \ell (m-1) + \ell_0,
$$
with $1 \leq \ell_0 \leq \ell$. Our assumption on $t$ implies that 
\begin{equation}\label{eq:parameterm}
m= o \left(\frac{\log n}{\log\log n }\right).
\end{equation}

From \eqref{eq:trace_identity}, to prove Theorem~\ref{lem:CLT}, it suffices to understand the moments of $X_{\ell,t}$ and the scalar term $z_{\ell,t}$. We start by calculating the moments of $X_{\ell,t}$. 
Expand $X_{\ell,t}$ as \[
 X_{\ell,t}
 =\sum_{x_1,\ldots,x_m\in[n]}
   Q^{(\ell)}(x_1,x_2)\cdots Q^{(\ell)}(x_{m-1},x_m)
   Q^{(\ell_0)}(x_m,x_1).
\] 
We set $\ell_j = \ell$ if $1 \leq j \leq m-1$ and $\ell_m = \ell_0$. Write $\mathbf x=(x_{ij})_{i\in[p],j\in[m]}$ and use the cyclic
convention $x_{i,m+1}=x_{i1}$. Then
\begin{equation}\label{eq:trace-entries}
 \E X_{\ell,t}^p
 =\sum_{\mathbf x\in[n]^{[p]\times[m]}}
   \E\prod_{i=1}^p\prod_{j=1}^m
      Q^{(\ell_j)}(x_{ij},x_{i,j+1}).
\end{equation}

From \eqref{eq:D}, we can further expand each $Q^{(\ell_j)}(x_{ij},x_{i,j+1})$ term as
\[
 Q^{(\ell_j)}(x_{ij},x_{i,j+1})
 =\underline P^{(\ell_j)}(x_{ij},x_{i,j+1})
  -\frac1N\sum_{k=1}^{\ell_j}R^{\ell_j,k}(x_{ij},x_{i,j+1}).
\]

A \emph{block} is one length-$\ell_j$ path representing the matrix
entry $Q^{(\ell_j)}(x_{ij},x_{i,j+1})$ at a fixed position $(i,j)$
in \eqref{eq:trace-entries}. The index $i$ selects a trace and
$j$ selects one of its $m$ consecutive blocks. Concatenating
these $m$ paths gives a closed path of length $t$.

Let $\mathcal I\subseteq[p]\times[m]$ be the set of positions
where a remainder contributing a factor $-1/N$ is chosen. For each $(i,j)\in\mathcal I$,
choose its omitted position $k_{ij}\in[\ell_j]$ (this determines a term involving $R^{\ell_j,k_{ij}}$), and write
$\mathbf k=(k_{ij})_{(i,j)\in\mathcal I}$. 

We call a block \emph{ordinary} if $(i,j)\notin\mathcal I$,
so its path comes from $\underline P^{(\ell_j)}$, and a
\emph{remainder block} if $(i,j)\in\mathcal I$, so its path
comes from $R^{\ell_j,k_{ij}}$.
We then obtain the following expansion:
\begin{align}\label{eq:block-choices}
 &\prod_{i=1}^p\prod_{j=1}^mQ^{(\ell_j)}(x_{ij},x_{i,j+1})\\
 =&\sum_{\mathcal I\subseteq[p]\times[m]}
       \sum_{\mathbf k\in\prod_{(i,j)\in\mathcal I}[\ell_j]}
       \left(\frac{-1}{N}\right)^{|\mathcal I|}
       \prod_{(i,j)\notin\mathcal I}
          \underline P^{(\ell_j)}(x_{ij},x_{i,j+1})
       \prod_{(i,j)\in\mathcal I}
          R^{\ell_j,k_{ij}}(x_{ij},x_{i,j+1}).
\end{align}

An \emph{ordered trace family} is the complete data
\[
 \gamma=\bigl((\gamma_{ij})_{i\in[p],j\in[m]},
                   \mathcal I,\mathbf k\bigr),\qquad
 \gamma_{ij}=(\mathbf e_{ij,s},\mathbf f_{ij,s})_{s=1}^{\ell_j},
\]
where $\mathcal I\subseteq[p]\times[m]$,
$\mathbf k\in\prod_{(i,j)\in\mathcal I}[\ell_j]$, and
\[
 \begin{aligned}
 \gamma_{ij}&\in\cT^{\ell_j}(x_{ij},x_{i,j+1})
                         &&(i,j)\notin\mathcal I,\\
 \gamma_{ij}&\in\cR^{\ell_j,k_{ij}}(x_{ij},x_{i,j+1})
                         &&(i,j)\in\mathcal I,
 \end{aligned}
\]
with  $x_{i,m+1}=x_{i1}$.
An ordered trace family corresponds to one term in the sum from \eqref{eq:block-choices}.
The \emph{omitted positions} are precisely the positions
$(i,j,k_{ij})$ for $(i,j)\in\mathcal I$.

The degree weight of the family $\gamma$ is the product of the
weights of its block paths $\gamma_{ij}$:
\begin{align}\label{eq:omega}
 \omega(\gamma)=\prod_{i=1}^p\prod_{j=1}^m\omega(\gamma_{ij})
 =\prod_{i=1}^p\prod_{j=1}^m\prod_{s=1}^{\ell_j}
       (d_{\mathbf e_{ij,s}}^+)^{-1}.
\end{align}
Its associated random matching factors have expectation
\begin{equation}\label{eq:joint-matching}
 \mu(\gamma)=\E\left[
  \prod_{(i,j)\notin\mathcal I}\prod_{s=1}^{\ell_j}
       \underline M_{\mathbf e_{ij,s}\mathbf f_{ij,s}}\prod_{(i,j)\in\mathcal I}
       \left(\prod_{s<k_{ij}}
          M_{\mathbf e_{ij,s}\mathbf f_{ij,s}}\right)
       \left(\prod_{s>k_{ij}}
          \underline M_{\mathbf e_{ij,s}\mathbf f_{ij,s}}\right)
 \right].
\end{equation}
Inserting all the definitions into \eqref{eq:block-choices}
and  \eqref{eq:trace-entries} gives the following exact expansion:
\begin{equation}\label{eq:trace-paths}
 \E X_{\ell,t}^p
 =\sum_{\mathcal I\subseteq[p]\times[m]}
     \sum_{\mathbf k\in\prod_{(i,j)\in\mathcal I}[\ell_j]}
     \left(\frac{-1}{N}\right)^{|\mathcal I|}
     \sum_{\mathbf x\in[n]^{[p]\times[m]}}
     \sum_{\substack{\gamma_{ij}\in
              \cT^{\ell_j}(x_{ij},x_{i,j+1})\\(i,j)\notin\mathcal I}}
     \sum_{\substack{\gamma_{ij}\in
              \cR^{\ell_j,k_{ij}}(x_{ij},x_{i,j+1})\\(i,j)\in\mathcal I}}
       \omega(\gamma)\mu(\gamma).
\end{equation}

\subsection{Grouping the terms in the trace sum}
We next separate the equality pattern of each path family from
its actual vertex and half-edge labels.

\begin{definition}[Canonical family]\label{def:canonical}
Explore the blocks $\gamma_{ij}$ in lexicographic order $(i,j)$,
reading the initial vertex and then the steps of each block
in increasing order.
Replace the distinct vertices by $1,2,\ldots,v$ in order of
first appearance. At each vertex, number the distinct heads
and the distinct tails separately starting from $1$, again in
order of first appearance. 

This relabeling includes the half-edges at omitted steps.
Keep the remainder-block set $\mathcal I$, the cuts
$\mathbf k$, and the order of all traces, blocks, and steps
unchanged. The resulting family is called the
\emph{canonical family}. We  denote it by
$\gamma$, with block paths $\gamma_{ij}$.
\end{definition}

Since there is one omitted edge in the definition of $R^{\ell,k}$ when we perform the path expansion, we will define the augmented support of a canonical family by putting back these edges.

\begin{definition}[Augmented support]\label{def:augmented-support}
Let $\gamma$ be a canonical family, with remainder blocks
$\mathcal I$ and cuts $\mathbf k$ as in \eqref{eq:trace-paths},
and put $r=|\mathcal I|$.
Its \emph{augmented support} is the directed multigraph with
all vertices of $\gamma$, one edge for each distinct oriented edge at a non-omitted step, and one formal directed edge for
each of the $r$ omitted steps.
Each formal edge is distinct from every other edge,
even when their endpoints coincide.

Let $a$ be the number of
distinct edges at non-omitted steps and $v$ the number of
vertices. Define
\begin{equation}\label{eq:q}
q=a+r-v.
\end{equation}
We call $q$ the \textit{excess} of the augmented support of $\gamma$. 
\end{definition}

\begin{lemma}\label{lem:excess}
     We have $a\geq v$. Hence $q\ge r$.
\end{lemma}
\begin{proof}
  Consider the support graph formed by all non-omitted edges
and all $v$ vertices. Each remainder block has genus at least
two by the definition of tangles, and the paths before and
after the omitted edge are tangle-free. There are two cases:
\begin{itemize}
    \item If the two paths are disjoint, adding the
    omitted edge back makes their union connected without
    creating a cycle. Then each path must contain at least
    one cycle for the remainder block to have genus at
    least two.
    \item If they intersect, their union is connected, and
    putting the omitted edge back increases its genus by
    at most one. Their union therefore already contains
    a cycle.
\end{itemize}
Any connected component containing an edge from a remainder
block therefore contains a cycle. Moreover, any connected
component containing no edges from remainder blocks contains
an entire closed trace and thus contains a cycle.
Thus every connected component in the support graph contains
a cycle; therefore $a\geq v$.
\end{proof}

Let $\deg x$ be the non-oriented degree of $x$ in the augmented support, with a
loop contributing two and parallel edges counted separately.
Each of its $a+r$ edges contributes two to the total degree, so
\begin{align}
 \sum_x\deg x=2(a+r),\qquad
 \sum_x(\deg x-2)
 =\sum_x\deg x-2v
 =2(a+r)-2v
 =2q. \label{eq:q_sum}
\end{align}

Let $a_1$ count consistent single edges whose unique
factor is centered, as in Lemma~\ref{lem:matching}.
For the canonical family $\gamma$, write
$(v_{ij,0},\ldots,v_{ij,\ell_j})$ for the sequence of vertex
symbols in block $\gamma_{ij}$.  For each $x\in[v]$, define the total number of path steps that leave vertex $x$
\[
 t_x=\#\bigl\{(i,j,s): i\in [p], j\in [m], s\in [\ell_j],
                         v_{ij,s-1}=x\bigr\}.
\] 
Since every step has one initial vertex, \[\sum_{x=1}^v t_x=pt.\]

For a fixed vertex tuple, grouping the factors
in \eqref{eq:omega} by the departure vertex therefore gives
\[
 \begin{aligned}
 \omega(\gamma)
 &=\prod_{i=1}^p\prod_{j=1}^m\prod_{s=1}^{\ell_j}
       (d_{\mathbf e_{ij,s}}^+)^{-1}
 =\prod_{x=1}^v(d_{i_x}^+)^{-t_x}.
 \end{aligned}
\]

Let $h_x^\pm$ count distinct heads and tails at vertex symbol
$x$, including those at omitted steps. Fix the first-visit tuple
$\mathbf i=(i_1,\ldots,i_v)$, so vertex symbol $x$ receives
the actual vertex label $i_x$. Its $h_x^+$ distinct head symbols must
receive distinct elements of $E^+(i_x)$. In their order of first
appearance, the first symbol has $d_{i_x}^+$ choices, the second
has $d_{i_x}^+-1$, and so on. Thus the number of assignments for heads is
\[
 d_{i_x}^+(d_{i_x}^+-1)\cdots(d_{i_x}^+-h_x^++1)
 =(d_{i_x}^+)_{h_x^+}.
\]
Similarly, assigning
the $h_x^-$ distinct tail symbols to $E^-(i_x)$ gives
$(d_{i_x}^-)_{h_x^-}$ choices. Multiplying
these counts therefore gives the total number of assignments:
\[
 \prod_{x=1}^v(d_{i_x}^+)_{h_x^+}(d_{i_x}^-)_{h_x^-}.
\]

For a canonical family $\gamma$, let $W(\gamma)$ denote the sum
of all terms in \eqref{eq:trace-paths} whose ordered path families
have canonical representative $\gamma$. We call $W(\gamma)$ the
\emph{class contribution}. This sum includes all choices of
distinct vertices $i_1,\ldots,i_v\in[n]$ and all half-edge
assignments; the block choices and cuts are fixed by $\gamma$.

Relabeling heads and tails preserves the uniform matching law,
so all these assignments have the same matching expectation
$\mu(\gamma)$. Multiplying the assignment count by the degree
weight and this expectation, including the coefficient $(-1/N)^r$,
and summing over the vertex tuple therefore gives
\begin{equation}\label{eq:class-weight}
 W(\gamma)=
 \left(\frac{-1}{N}\right)^r\mu(\gamma)
 \sum_{\substack{\mathbf i\in[n]^v\\\text{all distinct}}}
       \prod_{x=1}^v
       (d_{i_x}^+)_{h_x^+}(d_{i_x}^-)_{h_x^-}
                       (d_{i_x}^+)^{-t_x}.
\end{equation}

Write $\cF_{q,a_1}$ for the canonical families in
\eqref{eq:trace-paths} with a given $q$ and $a_1$. Grouping its terms gives the following identity:
\begin{equation}\label{eq:trace-classes}
 \E X_{\ell,t}^p
 =\sum_{q\ge0}\sum_{a_1\ge0}
       \sum_{\gamma\in\cF_{q,a_1}}W(\gamma).
\end{equation}
We first bound the classes with $q+a_1\ge1$,
then compute the leading class $q=a_1=0$.

\subsection{Path-counting estimates}
We adapt the analysis in
\cite{bordenave2015new,coste2021spectral}.
All degrees below refer to the augmented support.
A \textit{branch vertex} has degree at least three, or it is one of the starting vertices $x_{i,1}, 1\leq i\leq p$ of a trace.
A \emph{maximal non-branching path} is a directed sequence of
edges of positive length whose endpoints are branch vertices  
and whose internal vertices are not branch vertices.

\begin{lemma}[Branches and single edges]\label{lem:mixed-geometry}
For each canonical family with a given excess $q$, the following holds: 
\begin{enumerate}
    \item There are at most $2q+p$ branch vertices, with total degree
at most  $6q+2p$, and at most  $3q+p$ maximal non-branching paths.
\item There are at most $6q$ distinct inconsistent edges
and at most $q\ell$ raw single edges.
\end{enumerate}

\end{lemma}\label{lem:branch}
\begin{proof}
   (1)  From \eqref{eq:q_sum}, the number of vertices with degree greater than two is at most 
   \begin{align}\label{eq:deg_minus_2}
   \sum_{x} (\deg(x)-2)=2q.
   \end{align}
   Together with at most  $p$ degree-two starting vertices of a trace, there are at most $2q+p$ many branch vertices. Then the total degree satisfies
    \begin{align}\label{eq:branch_degree}
         \sum_{x \text{ branch vertices}} \deg(x) \leq \sum_{x} (\deg(x)-2) +2 \#\{ \text{branch vertices}\}\leq 6q+2p.
    \end{align}
    The endpoints of each maximal non-branching path are branch vertices. By counting degrees, there are at most $3q+p$ such paths.

(2) An inconsistent edge shares a head with another outgoing edge or a tail with another
incoming edge. This shared vertex has degree at least three. Therefore, from \eqref{eq:q_sum} and \eqref{eq:deg_minus_2}, the number of distinct inconsistent edges is bounded by \[\sum_{x:\deg(x)\geq 3} \deg(x)\leq \sum_{x:\deg(x)\geq 3} (\deg(x)-2) + 4q\leq  6q.\]
Each raw single edge belongs to a remainder block, so from Lemma~\ref{lem:excess}, there are at most $r\ell \leq q\ell$ such edges.
\end{proof}

At a vertex with one incoming and one outgoing edge, the paths use one head and one tail. If the vertex has label $i$ and the paths leave it $k$ times, its half-edge choices and departure weights contribute $d_i^+d_i^-(d_i^+)^{-k}$. 
Generalizing $\beta_3$ in \eqref{eq:defbeta3}, we define the following quantity to track its contribution 
\begin{equation}\label{eq:beta}
 \beta_k=\frac1N\sum_i\frac{d_i^-}{(d_i^+)^{k-1}},
 \qquad k\ge1.
\end{equation}

\begin{lemma}\label{lem:beta}
We have $\beta_1=1$, $\beta_2=\rho^2$, and for $k\geq 3$, 
\begin{equation}\label{eq:beta-bounds}
 \beta_k\le e^{-\theta} \rho^k,
\end{equation}
where $\theta$ is as in \eqref{eq:defbeta3}. Moreover, \eqref{eq:defbeta3} holds with $\theta = (\delta-\sqrt \delta)/\Delta $.
\end{lemma}
\begin{proof}
The first two identities follow from the definition. Since \eqref{eq:gap} implies that
$(d_i^+)^{-1}\le \delta^{-1} \leq \rho$, it gives for $k \geq 3$,
\[
 \beta_k\le\rho^{k-3}\frac1N\sum_i\frac{d_i^-}{(d_i^+)^2}
 = \rho^{k-3}\beta_3 \leq e^{-\theta} \rho^k.
\]

It remains to prove the last statement.  
For $x\in[1/\Delta,\rho]$, we have the affine bound $x^2\le(\rho+1/\Delta)x-\rho/ \Delta$.
We apply this to $x=1/d_i^+$ and average with weights $d_i^-/N$:
\[
 \beta_3\le \left(\rho+\frac 1 \Delta\right)\rho^2 -\frac{\rho} \Delta
 = \left(1-\frac{1/\rho^2-1/\rho}{\Delta}\right)\rho^3
       \le \left(1-\frac{\delta-\sqrt\delta}{\Delta}\right)\rho^3,
\]
where the last inequality follows from $\rho \leq 1/\sqrt\delta$. Since $(1-x) \leq e^{-x}$, we get \[\beta_3 \leq  \exp (-(\delta - \sqrt \delta )/ \Delta) \rho^3\] as claimed.
\end{proof}

\begin{lemma}[Class weight]\label{lem:mixed-weight}
If $pt\le\sqrt N$, then for the canonical family $\gamma$ with  given $a_1$ and $q$, we have
\begin{equation}\label{eq:mixed-weight}
 |W(\gamma)|
  \le C  \rho^{pt} \left( \frac{2^6\Delta^{10}}{n \rho^{\ell+6}} \right)^{q}   \left( \frac{3  pt}{\rho \sqrt{2n}}\right)^{a_1}.
\end{equation}
\end{lemma}
\begin{proof}
    Recall the definition of $\mu(\gamma)$ from \eqref{eq:joint-matching}. There are $pt -r$ matching factors, and $a$ distinct edges. In Lemma~\ref{lem:matching}, we have $b\leq 6q$ due to part (2) of Lemma~\ref{lem:mixed-geometry}. Applying Lemma~\ref{lem:matching}  to \eqref{eq:class-weight} gives
    \begin{align}
        |W(\gamma)| &\leq C2^{6q}\left( \frac{3  pt}{\sqrt N}\right)^{a_1}
        \frac{1}{N^{a+r}} \sum_{\substack{\mathbf i\in[n]^v\\\text{all distinct}}}
       \prod_{x=1}^v
       (d_{i_x}^+)_{h_x^+}(d_{i_x}^-)_{h_x^-}
                       (d_{i_x}^+)^{-t_x}.
    \end{align}

    Let $r_x$ be the in-degree of $x$ in the augmented support graph. Then 
    \[ \sum_{x} r_x=a+r.\]
    Dropping the requirement that the vertex labels be distinct gives
    \begin{align}
       \frac{1}{N^{a+r}} \sum_{\substack{\mathbf i\in[n]^v\\\text{all distinct}}}
       \prod_{x=1}^v
       (d_{i_x}^+)_{h_x^+}(d_{i_x}^-)_{h_x^-}
                       (d_{i_x}^+)^{-t_x}\leq \prod_{x=1}^v \frac{1}{N^{r_x}}\left( \sum_{i=1}^n (d_i^+)_{h_x^+} (d_i^-)_{h_x^{-}} (d_i^+)^{-t_x}\right)
    \end{align}
Denote \[\mu_x:=\frac{1}{N^{r_x}}\sum_{i=1}^n (d_i^+)_{h_x^+} (d_i^-)_{h_x^{-}} (d_i^+)^{-t_x}.\]
\begin{itemize}
    \item  Call a vertex \textit{ordinary} if it has degree two and has no incident formal edge; then $r_x=h_x^+=h_x^-=1$, and $\mu_x=\beta_{t_x}$. By Lemma~\ref{lem:beta}, if $t_x=1, \mu_x=1$ and $\mu_x\leq \rho^{t_x}$ if $t_x\geq 2$.  Then the number of ordinary vertices with $t_x=1$ is bounded by the total number of centered consistent single edges (which is $a_1$ by definition), centered inconsistent single edges (bounded by $6q$ from Lemma~\ref{lem:mixed-geometry}) and raw single edges (bounded by $\ell q$ from Lemma~\ref{lem:mixed-geometry}), which is  $a_1+(6+\ell)q$.
    
    \item All other vertices are called \textit{exceptional}, including vertices with degree at least three (whose total degree is at most $6q$ from Lemma~\ref{lem:mixed-geometry}) and degree-two vertices incident to a formal edge.  And for an exceptional $x$,
    \begin{align}
        \mu_x\leq  n^{1-r_x}\Delta^{h_x^++h_x^-}\rho^{t_x}\leq n^{1-r_x} \Delta^{\deg(x)} \rho^{t_x}.
    \end{align}
    Their total degree is at most $6q+4r$ by Lemma~\ref{lem:mixed-geometry}.
    Moreover, we have 
    \begin{align}
        \sum_{x~\text{exceptional}}1-r_x=v-(a+r)=-q.
    \end{align}
\end{itemize}
     Putting all estimates together, using $r\leq q$ from Lemma~\ref{lem:excess} and $N\geq \delta n \geq 2n$,
    \begin{align}
        |W(\gamma)|&\leq C 2^{6q}\Delta^{10q}n^{-q} \rho^{pt-a_1-(6+\ell)q}
 \left(\frac{3 pt}{\sqrt{2n}}\right)^{a_1}
    \end{align}
    as desired.
\end{proof}

Having bounded the weight of a single class, we now count how
many classes can have a given excess $q$.
\begin{lemma}[Canonical count]\label{lem:mixed-count}
\begin{equation}\label{eq:mixed-count} 
 \sum_{a_1\ge0}|\cF_{q,a_1}|\le(2pt)^{19pm(q+1)}.
\end{equation}
\end{lemma}

\begin{proof}
We adapt the spanning-tree and cycling-time encoding of
\cite{bordenave2015new,coste2021spectral}. Our support may be
disconnected, so we use a spanning forest. We record the support and the omitted edges first, then describe each tangle-free path by its crossings outside
the spanning forest and its repeated turns around a cycle. Finally, we
specify the half-edge identifications.

\noindent\emph{1. The support: at most $9q+3p$ integers.} For any $\gamma\in \mathcal F_{q,a_1}$, explore $\{\gamma_{ij}, i\in [p],j\in [m]\}$ in lexicographic order as in Definition~\ref{def:canonical}.
Number the vertices in the order of first visit from $1$ to $v$. Number the edges in the order of first visit from $1$ to $a+r$. 
Number the maximal non-branching paths in the order of the smallest edge label within the path.
By Lemma~\ref{lem:mixed-geometry}, there are at most $3q+p$ maximal
non-branching paths. Reserve $3q+p$ slots for triples
\[
 (\text{initial branch vertex},\ \text{final branch vertex},\ \text{length}).
\]
List the actual maximal non-branching paths in their assigned order and fill unused
slots with $(0,0,0)$.

\medskip
\noindent\emph{2. The block types, cuts, and formal edges: $2pm$ integers.} 
For each block, record
\[
 \begin{cases}
 (0,0),&\text{if it is ordinary},\\
 (k_{ij},\ \text{index assigned to the edge at the omitted step}),
     &\text{if it is a remainder}.
 \end{cases}
\]
The first entry distinguishes the block types and specifies the
cut when present; the second locates the formal edge. This uses
two integers per block, hence $2pm$ in total.

\medskip
\noindent\emph{3. Spanning trees and cycling times: at most $2pm(2q+5)$ integers.}
Start with all support vertices and no edges. Process the
support edges in their assigned order, ignoring orientations.
Add an edge if it does not create a cycle; otherwise, skip it.
The resulting spanning forest $T$ has one tree in each connected
component and is uniquely determined.
In a connected component $H$, exactly $g(H)$ edges lie outside $T$.
 Since
\[
 q=\sum_H\bigl(|E(H)|-|V(H)|\bigr)
   =\sum_H\bigl(g(H)-1\bigr),
\]
we have $g(H)\le q+1$ for every connected component in the augmented support graph.

For each ordinary block, consider its whole path.
For each remainder block $(i,j)$, consider the paths formed
by its first $k_{ij}-1$ and last $\ell_j-k_{ij}$ steps.
There are $pm+r\le2pm$ such indexed tangle-free paths,
ordered by block and then by their position in the block.

Denote one of these paths by $\gamma$.
Its length is determined by the block type and cut.
If $\gamma$ is not simple, tangle-freeness implies that its
support has exactly one unoriented cycle.
Deleting any edge outside this cycle disconnects the support, and we call such an edge a bridge.
After crossing a bridge, a directed walk cannot return to the connected
component it left, since this would require traversing the same edge
in the opposite direction. Thus
\[
 \gamma=\gamma_{\rm in}\mathbf c^\nu\mathbf c'\gamma_{\rm out},
 \qquad \nu\ge1,\quad 0\le|\mathbf c'|<|\mathbf c|,
\]
where $\gamma_{\rm in}$ and $\gamma_{\rm out}$ are simple,
possibly empty paths containing no cycle edges, $\mathbf c$ is
the directed cycle based at the entry vertex visited $\nu$ times, and $\mathbf c'$
is an initial segment of $\mathbf c$.
Delete the extra complete turns and set
\[
 \widehat\gamma
 =\gamma_{\rm in}\mathbf c\mathbf c'\gamma_{\rm out}.
\]
Every edge of $\widehat\gamma$ is traversed at most twice:
edges outside the cycle are traversed once, while cycle edges
are traversed once in $\mathbf c$ and possibly once more in
$\mathbf c'$.

For a simple or empty path, set $\widehat\gamma=\gamma$
instead. No additional cycle records are needed: the
reconstruction below restores the deleted turns from the
compressed path and the original length, which is fixed by
the block type and cut.

We now encode $\widehat\gamma$ using $T$. A
\emph{cycling time relative to $T$} is a time when the path
traverses an edge outside $T$.
Record the initial and final vertex indices of $\widehat\gamma$,
the number of times it traverses an edge outside $T$, and the
indices of these edges in traversal order. Record repeated
traversals of the same edge separately.
There are at most $q+1$ edges outside $T$ in its
component, and each is traversed at most twice. Thus the list has
at most $2(q+1)$ entries.

These records determine $\widehat\gamma$. Starting at its initial
vertex, follow the unique tree path to the initial vertex of the
first listed edge, traverse that edge, and continue in the same
way through the list. After the last listed edge, follow the unique
tree path to the recorded final vertex. If the list is empty,
follow the tree path between the two recorded vertices.
Thus one indexed path uses at most
\[
 \underbrace{2}_{\text{initial and final vertices}}
 +\underbrace{1}_{\text{list length}}
 +\underbrace{2(q+1)}_{\text{list entries}}
 =2q+5
\]
integers. The at most $2pm$ indexed paths therefore use at most
$2pm(2q+5)$ integers.

\medskip
\noindent\emph{4. The half-edge identifications: at most $6q$ integers.} The preceding records specify the support edges traversed by
each path, but not which distinct outgoing edges share a head
or which distinct incoming edges share a tail.

At each vertex of degree at least three, number the distinct
heads and tails separately, starting from one in order of first appearance. For each outgoing or incoming edge, record its head or tail number accordingly.
Equal numbers specify a shared half-edge. This uses one integer per outgoing edge and one per incoming
edge at a vertex of degree at least three. Thus the cost is at most
\[
 \sum_{\deg x\ge 3}\deg x\le 6q.
\]

\medskip
\noindent\emph{Reconstruction.}
First expand the nonzero support records into directed paths
and cycles. The block records identify the formal edges, and
the edge indices in the support determine $T$. For each indexed path, its
endpoints and list of edges traversed outside $T$ recover
$\widehat\gamma$ by the unique tree paths.

If $\widehat\gamma$ is simple, no turns were deleted, because
the compression keeps one complete turn whenever a cycle is
present. Otherwise, its first repeated vertex identifies the
cycle $\mathbf c$ and its entry point. The original length
$|\gamma|$ is known from the block data, so insert
 $\frac{|\gamma|-|\widehat\gamma|}{|\mathbf c|}$
additional complete turns at that entry point to recover
$\gamma$.  Join the two paths of each remainder block through
its formal edge, and apply the head and tail identifications
from Step~4. Thus every
valid record list determines a unique canonical family.

\medskip
\noindent\emph{Counting the records.}
The total number of integers is at most
 \begin{align}
 &(9q+3p)+2pm+2pm(2q+5)+6q
 \le 19pm(q+1),
 \end{align}

The support has at most $p t$ edges and vertices, so every
vertex, edge, and half-edge identification index is at most
$pt$. Cuts and the lengths of the lists of edges traversed
outside $T$ are also at most $p t$. Including zero for unused
records, each integer therefore has at most $ pt+1\leq 2pt$ possible
values. 
Since distinct canonical families have distinct integer sequences,
this proves \eqref{eq:mixed-count}.
\end{proof}

\begin{figure}[ht]
\centering
\small
\vspace{-1ex}\par\noindent
\begin{minipage}[t]{.48\textwidth}
\vspace{0pt}
\centering
\begin{tikzpicture}[x=.9cm,y=.9cm,>=Stealth,
 vertex/.style={circle,fill=black,inner sep=1.6pt},
 start/.style={circle,draw,fill=white,inner sep=2pt},
 tree/.style={->,thick},
 outside/.style={->,thick,FireBrick},
 every node/.style={font=\small}]
 \node[start,label=above:$1$] (v1) at (0,0) {};
 \node[start,label=above:$2$] (v2) at (2,0) {};
 \node[vertex,label=above:$3$] (v3) at (4,0) {};
 \node[vertex,label=above:$4$] (v4) at (5.4,1) {};
 \node[vertex,label=right:$5$] (v5) at (6.8,0) {};
 \node[vertex,label=below:$6$] (v6) at (5.4,-1) {};
 \draw[tree] (v1) to[bend left=32] node[above] {$e_1$} (v2);
 \draw[outside] (v2) to[bend left=32] node[below] {$e_8$} (v1);
 \draw[tree] (v2) to[bend left=32] node[above] {$e_2$} (v3);
 \draw[outside] (v3) to[bend left=32] node[below] {$e_7$} (v2);
 \draw[->,thick,RoyalBlue] (v3) to[bend left=95]
   node[below] {$e_9$} (v2);
 \draw[tree] (v3) to[out=80,in=180] node[above left] {$e_3$} (v4);
 \draw[tree] (v4) to[out=0,in=100] node[above right] {$e_4$} (v5);
 \draw[tree] (v5) to[out=-100,in=0] node[below right] {$e_5$} (v6);
 \draw[outside] (v6) to[out=180,in=-80] node[below left] {$e_6$} (v3);
 \pgfresetboundingbox
 \path[use as bounding box] (-.4,-1.6) rectangle (7.4,1.6);
\end{tikzpicture}\par\smallskip
\makebox[\linewidth][l]{\textbf{\strut(a) Augmented support}}
\end{minipage}\hfill
\begin{minipage}[t]{.48\textwidth}
\vspace{0pt}
\centering
\begin{tikzpicture}[x=1cm,y=1cm,>=Stealth,
 vertex/.style={circle,fill=black,inner sep=1.6pt},
 every node/.style={font=\small},every path/.style={thick}]
 \node[vertex,label=above:$1$] (b1) at (0,0) {};
 \node[vertex,label=above:$2$] (b2) at (1.5,0) {};
 \node[vertex,label=above:$3$] (b3) at (3,0) {};
 \draw[->] (b1) to[bend left=35] node[above] {$1$} (b2);
 \draw[->] (b2) to[bend left=35] node[below] {$1$} (b1);
 \draw[->] (b2) to[bend left=35] node[above] {$1$} (b3);
 \draw[->] (b3) to[bend left=35] node[below] {$1$} (b2);
 \draw[->,RoyalBlue] (b3) .. controls (3,-1.2) and (1.5,-1.2)
   .. node[below,pos=.5] {$1$} (b2);
 \draw[->] (b3) .. controls (4.5,1.5) and (4.5,-1.5)
   .. node[right,pos=.5] {$4$} (b3);
 \pgfresetboundingbox
 \path[use as bounding box] (-1.3,-1.44) rectangle (5.72,1.44);
\end{tikzpicture}\par\smallskip
\makebox[\linewidth][l]{\textbf{\strut(b) Maximal non-branching paths}}
\end{minipage}

\par\medskip\noindent
\begin{minipage}[t]{.48\textwidth}
\vspace{0pt}
\centering
\begin{tikzpicture}[x=.95cm,y=.85cm,>=Stealth,
 vertex/.style={circle,draw,minimum size=7mm,inner sep=1pt},
 half/.style={rectangle,draw,fill=white,minimum size=3.5mm,inner sep=.5pt,font=\scriptsize},
 every node/.style={font=\small}]
 \node[vertex] (v) at (0,0) {$2$};
 \node[half] (t1) at (-1,.6) {$1$};
 \node[half] (t2) at (-1,-.6) {$2$};
 \node[half,draw=FireBrick,text=FireBrick,fill=white] (h1) at (1,0) {$1$};
 \draw[->,thick] (t1)--(v);
 \draw[->,thick] (t2)--(v);
 \draw[->,very thick,FireBrick] (v)--(h1);
 \draw[->] (-2.3,1) to[out=0,in=160] (t1);
 \draw[->] (-2.3,-1) to[out=0,in=-160] (t2);
 \draw[->,RoyalBlue] (-2.3,0) to[out=0,in=160] (t2);
 \draw[->,FireBrick] (h1) to[out=25,in=180] (2.3,.9);
 \draw[->,FireBrick] (h1) to[out=-25,in=180] (2.3,-.9);
 \node[above] at (-2.2,1) {$e_1$};
 \node[below] at (-2.2,-1) {$e_7$};
 \node[above,RoyalBlue] at (-2.2,0) {$e_9$};
 \node[above] at (2.2,.9) {$e_2$};
 \node[below] at (2.2,-.9) {$e_8$};
 \pgfresetboundingbox
 \path[use as bounding box] (-3.7,-2) rectangle (3.7,2);
\end{tikzpicture}\par\smallskip
\makebox[\linewidth][l]{\textbf{\strut(c) A shared head at vertex $2$}}
\end{minipage}\hfill
\begin{minipage}[t]{.48\textwidth}
\vspace{0pt}
\centering
\begin{tikzpicture}[x=1cm,y=1cm,>=Stealth,
 visit/.style={circle,draw,fill=white,minimum size=5mm,inner sep=1pt},
 every node/.style={font=\small}]
 \foreach \i/\v in {0/2,1/1,2/2,3/1,4/2}
   {\node[visit] (g\i) at ({.8*\i},.95) {$\v$};}
 \node[left] at (-.45,.95) {$\gamma_{22}$};
 \draw[->,thick,FireBrick] (g0)--node[above] {$e_8$} (g1);
 \draw[->,thick] (g1)--node[above] {$e_1$} (g2);
 \draw[->,thick,FireBrick] (g2)--node[above] {$e_8$} (g3);
 \draw[->,thick] (g3)--node[above] {$e_1$} (g4);
 \draw[->,thick] (1.6,.5)--node[right,align=left,font=\scriptsize]
   {} (1.6,-.45);
 \foreach \i/\v in {0/2,1/1,2/2}
   {\node[visit] (h\i) at ({.8+.8*\i},-.95) {$\v$};}
 \node[left] at (-.45,-.95) {$\widehat\gamma_{22}$};
 \draw[->,thick,FireBrick] (h0)--node[above] {$e_8$} (h1);
 \draw[->,thick] (h1)--node[above] {$e_1$} (h2);
 \pgfresetboundingbox
 \path[use as bounding box] (-2.3,-1.7) rectangle (4.73,1.7);
\end{tikzpicture}\par\smallskip
\makebox[\linewidth][l]{\textbf{\strut(d) Cycle compression in $\gamma_{22}$}}
\end{minipage}

\caption{Encoding two traces with $p=2$, $m=2$, $\ell=4$, and $t=8$:
$(1,2,3,4,5)\mid(5,6,3,2,1)$ and
$(2,3,2,1,2)\mid(2,1,2,1,2)$, where bars separate blocks.
The only remainder block is $\gamma_{21}$, with $k_{21}=2$.
The matching factor at $3\to2$ is omitted, leaving the tangle-free
paths $2\to3$ and $2\to1\to2$.
In (a), the blue  edge $e_9$ records this omitted step;
it is distinct from $e_7$ but uses the same head and tail.
Thus $a=8$, $r=1$, $v=6$, and $q=3$.
Open circles are the trace starting vertices; the red edges
and $e_9$ lie outside the spanning tree $T=\{e_1,e_2,e_3,e_4,e_5\}$.
In (b), each arrow represents a maximal non-branching path
labelled by its length. Panel (c) shows the  half-edge
labels at vertex $2$, with heads and tails numbered separately:  $e_2,e_8$ share a head and $e_7,e_9$ share a tail.
In (d), deleting one complete turn from
$\gamma_{22}=(2,1,2,1,2)$ gives $\widehat\gamma_{22}=(2,1,2)$.}
\label{fig:canonical-encoding}
\end{figure}
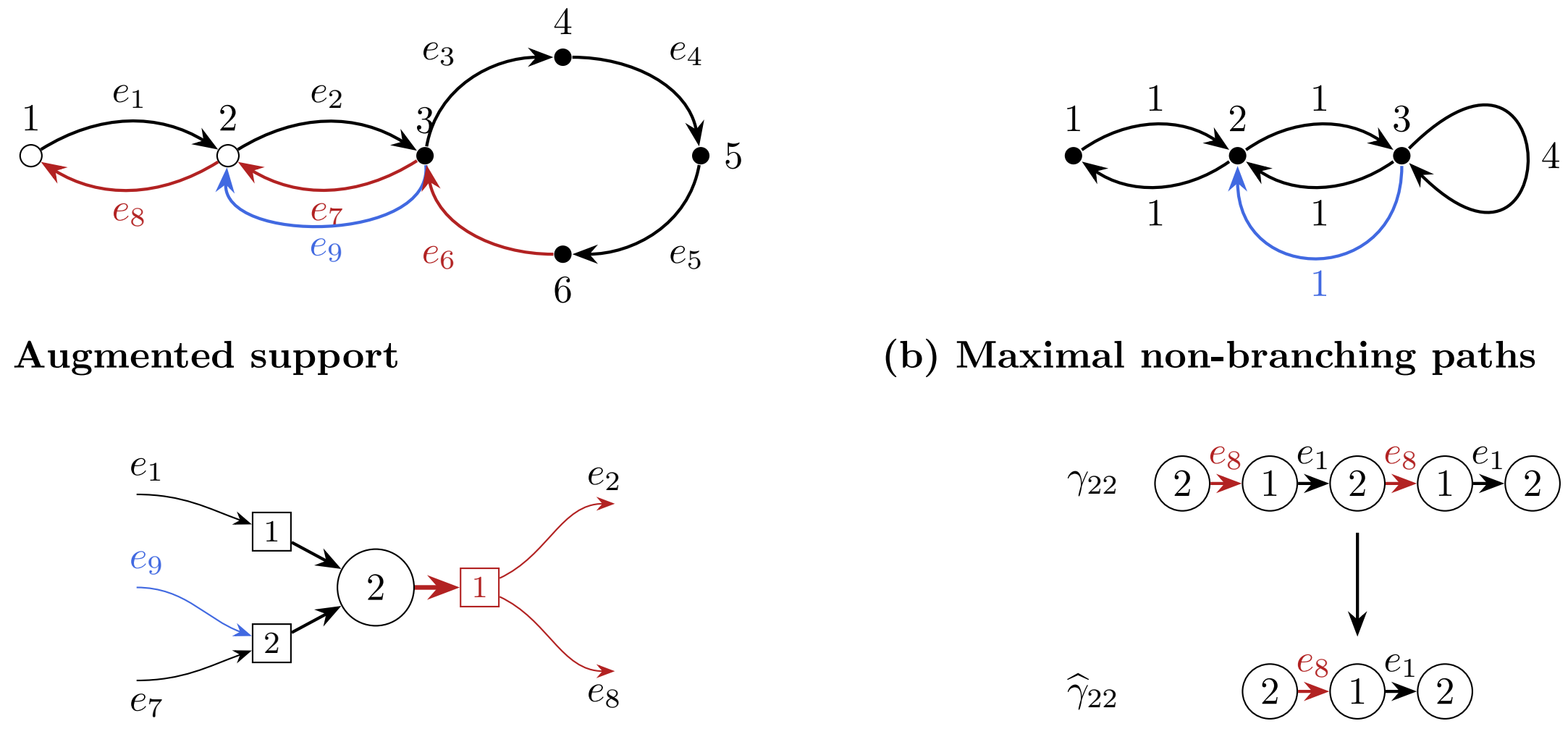

Figure~\ref{fig:canonical-encoding} illustrates the encoding for two
traces according to  Lemma~\ref{lem:mixed-count}.

We combine Lemmas~\ref{lem:mixed-weight} and~\ref{lem:mixed-count}
to control the classes with $q+a_1\ge1$. 
\begin{lemma}[Error sum]\label{lem:mixed-errors}
For $t,\ell$ in \eqref{eq:parametert}-\eqref{eq:parameters}
and any fixed positive integer $p$, the following holds:
\begin{equation}\label{eq:error-sum}
 \sum_{\substack{q,a_1\ge0\\q+a_1\ge1}}
     \sum_{\gamma\in\cF_{q,a_1}}|W(\gamma)|
 \le\frac{24Cpt\sqrt{\Delta}\rho^{ p t}}{\sqrt n}
                 (2 p t )^{19pm}
 =o(\rho^{p t}).
\end{equation}
\end{lemma}
\begin{proof}
Since $N\ge2n$, the parameter choice \eqref{eq:parametert} gives
\[
 p t =O_p((\log n)^2)=o(\sqrt N),
 \quad \hbox{ and } \quad \Delta^\ell\le n^{1/8}.
\]
Also, since $m = o (\log n /  \log \log n)$ from \eqref{eq:parameterm}, we find
\begin{align}\label{eq:v}
 pm\log(2 pt)
 =o(\log n).
\end{align}
Thus $(2 pt)^{pm}=n^{o(1)}$.
By \eqref{eq:gap}, $\rho^2\ge\Delta^{-1}$. Together with
$\Delta^\ell\le n^{1/8}$, this gives
\[
 \frac{2^6\Delta^{10}}{n\rho^{\ell+6}}
 \le\frac{2^6\Delta^{13}}{n^{15/16}}.
\]
Multiplying Lemmas~\ref{lem:mixed-weight} and ~\ref{lem:mixed-count} 
therefore yields, for every $q,a_1\ge0$,
 \begin{align}
 \rho^{- pt}\sum_{\gamma\in\cF_{q,a_1}}|W(\gamma)|
 &\le C(2pt)^{19pm(q+1)}
 \left(\frac{2^6\Delta^{13}}{n^{15/16}}\right)^q
 \left(\frac{3pt\sqrt\Delta}{\sqrt{2n}}\right)^{a_1}. \label{eq:normalized_qa}
 \end{align}
Set
\[
  u=\frac{2^6\Delta^{13}(2pt)^{19pm}}{n^{15/16}},
 \qquad
 v=\frac{3pt\sqrt\Delta}{\sqrt{2n}}.
\]
Then from \eqref{eq:v} and the condition $\Delta=n^{o(1)}$, for sufficiently large $n$, \[0<u\le v=n^{-1/2+o(1)}\le1/2,\] 
and the bound for fixed $(q,a_1)$ for the right hand side of \eqref{eq:normalized_qa} is
$C(2pt)^{19pm}u^qv^{a_1}$.
By the geometric sum,
\[
 \begin{aligned}
 \sum_{\substack{q,a_1\ge0\\q+a_1\ge1}}u^qv^{a_1}
 =\left(\sum_{q\ge0}u^q\right)
   \left(\sum_{a_1\ge0}v^{a_1}\right)-1=\frac1{(1-u)(1-v)}-1=\frac{u+v-uv}{(1-u)(1-v)}
 \le8v.
 \end{aligned}
\]
Therefore the left hand side of \eqref{eq:normalized_qa} satisfies
\[
 \begin{aligned}
 \rho^{- pt}
 \sum_{\substack{q,a_1\ge0\\q+a_1\ge1}}
       \sum_{\gamma\in\cF_{q,a_1}}|W(\gamma)|
 &\le \frac{24Cpt\sqrt{\Delta}(2pt)^{19pm}}{\sqrt n}=o(1).
 \end{aligned}
\]
 This proves \eqref{eq:error-sum}.
\end{proof}

\subsection{The leading term}
Lemma~\ref{lem:mixed-errors} and \eqref{eq:trace-classes} give,
for every fixed positive integer $p$,
\begin{equation}\label{eq:leading-reduction}
 \E X_{\ell,t}^p
 =\sum_{\gamma\in\cF_{0,0}}W(\gamma)
       +o(\rho^{pt}).
\end{equation}

When $q=0$, every  vertex in the canonical family $\gamma$ has in-degree and out-degree
one in the augmented support, so each connected component is a directed cycle with $r=0$ and all matching factors are centered. All edges are consistent, because each vertex has only
one incoming and one outgoing edge. Since $a_1=0$, every
edge has multiplicity at least two.

Each trace lies in one cycle component. Let $\mathcal Q$ be
the resulting partition of $[p]$: a set $A\in\mathcal Q$
contains the indices of the traces following the same cycle.
If this cycle has length $j_A$, then
\[
 j_A\mid t ,\qquad
 t_A=\frac{t|A|}{j_A}\ge2,\qquad
 a=\sum_{A\in\mathcal Q}j_A\le\frac{pt}{2}.
\]
Here $t_A$ is the multiplicity of every edge in that cycle,
and $a$ is the total number of distinct edges. 

For all $\gamma\in \mathcal F_{0,0}$, with a given
$\mathcal Q,\mathbf j=(j_A)_{A\in\mathcal Q}$,
Lemma~\ref{lem:matching} gives for some  $\varepsilon_{\mathcal Q,\mathbf{j}}$ depending  on $\mathcal Q, \mathbf{j}$, 
\begin{align}
    \mu(\gamma)=\frac{1+\varepsilon_{\mathcal Q,\mathbf{j}}}{(N)_a}, \qquad |\varepsilon_{\mathcal Q,\mathbf{j}}| \leq \frac{2pt}{N}.
\end{align}

Consider the traces following the same cycles in $A$. Since their starting vertices can be different, there are  $j_A^{|A|-1}$ choices for their relative starting locations. List the actual vertices of this cycle in directed order, starting at the
initial vertex of the trace corresponding to the minimal index in $A$:
\[
 i_{A,1}\to\cdots\to i_{A,j_A}\to i_{A,1}.
\]
choose  distinct values of $i_{A,u}, A\in \mathcal Q, 1\leq u\leq j_A$ determines  the labeling of all $\gamma\in \mathcal F_{0,0}$.
Then we have the following full expansion:
 \begin{align}\label{eq:joint-cycles}
 &\sum_{\gamma\in\cF_{0,0}}W(\gamma)\\
 =&\sum_{\mathcal Q\ {\rm partition\ of}\ [p]}
     \sum_{\substack{\mathbf j=(j_A)_{A\in\mathcal Q}\\
                       j_A\mid t,\  t |A|  /j_A\ge2}}
       \frac{1+\varepsilon_{\mathcal Q,\mathbf j}}{(N)_a}
       \prod_{A\in\mathcal Q}j_A^{|A|-1}\sum_{\substack{(i_{A,u})\in[n]^a\\\text{all distinct}}}
       \prod_{A\in\mathcal Q}\prod_{u=1}^{j_A}
       \frac{d_{i_{A,u}}^-}{(d_{i_{A,u}}^+)^{t_A-1}}.
 \end{align}

We now evaluate the pairings of traces on cycles of length
$t$ and show that all other cycle arrangements are negligible.

\begin{lemma}[Gaussian moments]\label{lem:moments} For $t,\ell$ as in \eqref{eq:parametert}-\eqref{eq:parameters}, $\Delta = n^{o(1)}$,
and any fixed positive integer $p$,
\begin{equation}\label{eq:moments}
 \E\left(\frac{X_{\ell,t}}{\sqrt{t}\rho^{t}}\right)^p
 \longrightarrow
 \begin{cases}
 (p-1)!!,&p\text{ even},\\
 0,&p\text{ odd}.
 \end{cases}
\end{equation}
Therefore \[
 \frac{X_{\ell,t}}{\sqrt{t}\,\rho^{t}}
 \Longrightarrow\mathcal N(0,1).
\]
\end{lemma}
\begin{proof}
Since $p$ is fixed and $t=O((\log n)^2)$, we have $pt=o(\sqrt N)$.
From \eqref{eq:Na_ratio}, uniformly for $a\le pt/2$, 
\begin{align}\label{eq:factorial-bound}
\frac{N^a}{(N)_a}=1+o(1).
\end{align}

\medskip
\noindent\emph{Leading terms.}
For even $p$, choose a pairing $\mathcal Q$, that is, a partition
of $[p]$ into $p/2$ sets of size two. Set $j_A=t$ for every
pair $A\in\mathcal Q$. The two traces in each pair traverse their
common cycle once, so every edge is used twice:
\[
 t_A=\frac{t|A|}{j_A}
     =\frac{2t}{t}=2,
 \qquad
 a=\sum_{A\in\mathcal Q}j_A=\frac{pt}{2}.
\]
All pairings with each trace traversing a cycle of length $t$ give the leading term in \eqref{eq:joint-cycles}.

Enumerate the labels of these $a$ cycle vertices as
$i_1,\ldots,i_a$. Since $t_A-1=1$, the degree sum in
\eqref{eq:joint-cycles} for this pairing is
\[
 \sum_{\substack{(i_1,\ldots,i_a)\in[n]^a\\\text{all distinct}}}
       \prod_{u=1}^a\frac{d_{i_u}^-}{d_{i_u}^+}.
\]
If we drop the distinct-label restriction, each label ranges
independently over $[n]$, and the sum factors:
\[
 \sum_{i_1=1}^n\cdots\sum_{i_a=1}^n
       \prod_{u=1}^a\frac{d_{i_u}^-}{d_{i_u}^+}
 =\prod_{u=1}^a\left(\sum_{i=1}^n\frac{d_i^-}{d_i^+}\right)=(N\rho^2)^a.
\]
We now show that imposing distinctness changes this sum by a
relative $o(1)$.  Assumption~\eqref{eq:gap} gives
\[
 \max_i\frac{d_i^-}{d_i^+N\rho^2}
 \le\frac{\Delta \delta}{N}\le\frac{\Delta}{n}.
\]
Every tuple excluded by the distinct-label restriction satisfies
$i_s=i_{s'}$ for some $1\le s<s'\le a$. Since all summands are
nonnegative,
\[
 \begin{aligned}
 0&\le (N\rho^2)^a-
 \sum_{\substack{(i_1,\ldots,i_a)\in[n]^a\\\text{all distinct}}}
       \prod_{u=1}^a\frac{d_{i_u}^-}{d_{i_u}^+}\le\sum_{1\le s<s'\le a}
 \sum_{\substack{(i_1,\ldots,i_a)\in[n]^a\\i_s=i_{s'}}}
       \prod_{u=1}^a\frac{d_{i_u}^-}{d_{i_u}^+}.
 \end{aligned}
\]
For each pair $s<s'$, the two equal labels contribute a squared
factor, while the other $a-2$ labels can be summed separately:
\[
 \sum_{\substack{(i_1,\ldots,i_a)\in[n]^a\\i_s=i_{s'}}}
       \prod_{u=1}^a\frac{d_{i_u}^-}{d_{i_u}^+}
 =(N\rho^2)^{a-2}
       \sum_{i=1}^n\left(\frac{d_i^-}{d_i^+}\right)^2\le(N\rho^2)^{a-1}\max_i\frac{d_i^-}{d_i^+}\le\frac{\Delta}{n}(N\rho^2)^a.
\]
Summing over the $\binom a2$ pairs  gives
\[
 0\le1-
 \frac{\displaystyle
       \sum_{\substack{(i_1,\ldots,i_a)\in[n]^a\\\text{all distinct}}}
                         \prod_{u=1}^a\frac{d_{i_u}^-}{d_{i_u}^+}}
      {(N\rho^2)^a}
 \le\binom a2\frac{\Delta}{n}=o(1).
\]
By \eqref{eq:factorial-bound} and the matching error bound,
the contribution of this fixed pairing is
\[
 \frac{1+\varepsilon_{\mathcal Q,\mathbf j}}{(N)_a}
       t^{p/2}
       \sum_{\substack{(i_1,\ldots,i_a)\in[n]^a\\\text{all distinct}}}
                          \prod_{u=1}^a\frac{d_{i_u}^-}{d_{i_u}^+}
 =t^{p/2}\rho^{pt}(1+o(1)).
\]
There are $(p-1)!!$ such pairings when $p$ is even, and zero
when $p$ is odd.

\medskip
\noindent\emph{All other terms.}
For any given $\mathcal Q$, dropping label distinctness in a term of \eqref{eq:joint-cycles}
factors its degree sum as
\[
 \prod_{A\in\mathcal Q}
       \left(\sum_i\frac{d_i^-}{(d_i^+)^{t_A-1}}\right)^{j_A}
 =N^a\prod_{A\in\mathcal Q}\beta_{t_A}^{j_A}.
\]
Together with \eqref{eq:factorial-bound}, this bounds the
absolute value of the contribution for a given $\mathcal Q$ in  \eqref{eq:joint-cycles} by
\begin{equation}\label{eq:cycle-upper}
 \sum_{\substack{\mathbf j=(j_A)_{A\in\mathcal Q}\\
                       j_A\mid t,\ |A|t/j_A\ge2}} 2\prod_{A\in\mathcal Q}j_A^{|A|-1}\beta_{t_A}^{j_A}
\end{equation}
for sufficiently large $n$.

In \eqref{eq:moments}, $\mathbb EX_{\ell,t}^p$ is normalized by
$t^{p/2}\rho^{pt}$. Since
$\sum_{A\in\mathcal Q}|A|=p$, this normalization factors as
\[
 t^{p/2}\rho^{pt}
 =\prod_{A\in\mathcal Q}
       t^{|A|/2}\rho^{t|A|}.
\]
We may therefore divide the factor for a cycle followed by $b$
traces by $t^{b/2}\rho^{bt}$, and then multiply
these normalized factors over all cycles in \eqref{eq:joint-cycles}.  

We now separate the remaining terms into two cases.

\medskip
\noindent\emph{(i) $\mathcal Q$ is a pairing, but at least one
$j_A<t$.}
Here $p$ is even and every cycle is followed by two traces.
For the shorter cycles, define their normalized contribution by
\[
 s_2=\frac1{t\,\rho^{2t}}
       \sum_{\substack{j\mid t \\j<t}}
                    j\beta_{2t/j}^{j}.
\]
A proper divisor $j$ of $t$ satisfies $j\le t/2$,
so the multiplicity $2t/j$ is at least four.
Lemma~\ref{lem:beta} gives
\[
 j\beta_{2t/j}^{j}
\le\rho^{2t}j e^{-\theta j}
 \le\theta^{-1} \rho^{2t} ,
\]
where  the first inequality uses $1-x \leq e^{-x}$ for $x\geq 0$ and the last inequality uses $e^x\ge x$ for $x\geq 0$.

Write $\tau(t)$ for the number of positive divisors of $t$.
We use the standard bound \cite[Theorem 315]{hardy1979introduction} that  $\tau(t)=t^{o(1)}$. Hence, since $t \theta \geq t^{1/2 -o(1)}$ from \eqref{eq:parametert}, we get
\[
 0\le s_2\le
 \frac{\tau(t)}{t\theta}=  o(1).
\]
For one pair of traces, separate the length $j=t$ from
the shorter lengths. Since $\beta_2=\rho^2$, the sum of its
normalized cycle bounds is
\[
 \sum_{j\mid t}
       \frac{j\beta_{2t/j}^{j}}{t\,\rho^{2t}}
 =\frac{t\beta_2^{t}}
          {t\rho^{2t}}+s_2=1+s_2.
\]

Now fix a pairing $\mathcal Q$, which has $p/2$ pairs.
Summing the product of these bounds over the cycle length of
each pair gives
\[
\sum_{\substack{(j_A)_{A\in\mathcal Q}\\
                  j_A\mid t}}
       \prod_{A\in\mathcal Q}
       \frac{j_A\beta_{2t/j_A}^{j_A}}
            {t\rho^{2t}}=
       \prod_{A\in\mathcal Q}
       \left(\sum_{j\mid t}
          \frac{j\beta_{2t/j}^{j}}
               {t\rho^{2t}}\right)
       =(1+s_2)^{p/2}.
\]
The unique choice $j_A=t$ for every pair contributes
$1$ to this sum of bounds. Subtracting it leaves precisely
the choices with at least one shorter cycle:
\[
 (1+s_2)^{p/2}-1
 =o(1).
\]
Finally, \eqref{eq:cycle-upper} gives the factor $2$,
and there are $(p-1)!!$ perfect pairings. Their total normalized
absolute contribution with at least one shorter cycle is
therefore at most
\[
 2(p-1)!!\bigl[(1+s_2)^{p/2}-1\bigr]=o(1).
\]

\medskip
\noindent\emph{(ii) $\mathcal Q$ is not a pairing.}
Such a partition has at least one set $A$ with $|A|\ne2$.
For a cycle followed by $b\ne2$ traces, define the normalized
sum over all its permitted lengths by
\[
 s_b=\frac1{t^{b/2}\rho^{bt}}
       \sum_{\substack{j\mid t\\bt/j\ge2}}
                    j^{b-1}\beta_{bt/j}^{j}.
\]
For $b=1$, Lemma~\ref{lem:beta} immediately gives
\[
 0\le s_1\le\frac{\tau(t)}{\sqrt{t}}=o(1).
\]
For $b\ge3$, every permitted multiplicity is at least three and Lemma~\ref{lem:beta} gives
\[
 j^{b-1}\beta_{bt/j}^{j}
 \le\rho^{bt}j^{b-1}e^{-\theta j}
 \le\rho^{bt}(b-1)! \theta^{-(b-1)},
\]
where the last inequality uses $e^x\ge x^{b-1}/(b-1)!$.
Summing over the permitted lengths yields, for each fixed $b\ge3$,
\[
 0\le s_b\le
\frac{\tau(t)(b-1)!}{\sqrt t} ( \sqrt t \theta)^{-(b-1)}= o(1),
\]
where we have used $b \leq p$ and $(\sqrt t \theta)^{-1} \le t^{o(1)}$.
Cycles followed by two traces still have normalized length
sum $1+s_2$. Thus \eqref{eq:cycle-upper} bounds the total
normalized absolute contribution of the partitions in this case by
\begin{equation}\label{eq:nonleading-cycles}
 2\sum_{\substack{\mathcal Q\ {\rm partition\ of}\ [p]\\
                  \mathcal Q\ {\rm not\ a\ pairing}}}
 \left(\prod_{\substack{A\in\mathcal Q\\|A|=2}}(1+s_2)\right)
 \left(\prod_{\substack{A\in\mathcal Q\\|A|\ne2}}s_{|A|}\right)
 =o(1).
\end{equation}
Combining this bound with the leading pairing contribution gives
\[
 \sum_{\gamma\in\cF_{0,0}}W(\gamma)
 =t^{p/2}\rho^{pt}
 \begin{cases}
 (p-1)!!+o(1),&p\text{ even},\\
 o(1),&p\text{ odd}.
 \end{cases}
\]
Together with \eqref{eq:leading-reduction}, this proves
\eqref{eq:moments}: its error after normalization is
$o(t^{-p/2})=o(1)$.
Since the Gaussian distribution is determined by its moments,
the convergence in distribution holds.
\end{proof}

\subsection{Bounding the scalar product term}
To use \eqref{eq:trace_identity}, we must show that the
scalar product correction $z_{\ell,t}$ is negligible at the desired scale. It suffices to bound the second moment of $z_k$ for $1\leq k\leq \ell$.

\begin{lemma}[Negligible scalar term]\label{lem:scalar}
Assume $\Delta =n^{o(1)}$. For any integer sequence $k =k_n$ with $1 \leq k \leq \ell$ and $\ell$ as in \eqref{eq:parameters},
\begin{equation}\label{eq:scalar}
 \frac{z_k}{\rho^{k}}\longrightarrow0
       \quad\text{in probability}.
\end{equation}
\end{lemma}
\begin{proof}
\noindent\emph{Expansion and the support.}
Expand the two matrix entries into paths and group by canonical
families, as in \eqref{eq:trace-classes}:
\begin{equation}\label{eq:scalar-paths}
 \E z_k^2
 =\frac1{n^2}\sum_{x_1,y_1,x_2,y_2\in[n]}
       \E Q^{(k)}(x_1,y_1)Q^{(k)}(x_2,y_2)
 =\sum_\gamma W(\gamma).
\end{equation}
Each family $\gamma$ now consists of two  length-$k$
paths $\gamma_1,\gamma_2$, and
$W(\gamma)$ is given by \eqref{eq:class-weight} with the additional
factor $n^{-2}$. 

For this proof, a branch vertex is a vertex of degree at
least three or one of the endpoints $x_1,y_1,x_2,y_2$.
Maximal non-branching paths are defined as before. Moreover, $a$ counts distinct edges carrying
matching factors, $r\le2$ counts omitted edges,
$v$ counts vertices, and $a_1$ counts consistent centered single
edges. The augmented support has $a+r$ edges and $c$
connected components for $1\leq c\leq 2$. Define
\[
 q=a+r-v+2=\sum_{H\ {\rm component}}g(H)+2-c.
\]

By the same argument as in Lemma~\ref{lem:excess}, every connected component with edges from the remainder block contains a cycle. Then the number of tree components is at most $(2-r)$ and $a-v\geq -(2-r)=r-2$. This implies $q\geq 2r$.

If $q=0$, then $c=2$ and both components are trees. A directed
walk in a tree cannot repeat a vertex, and neither block can
be a tangled remainder block. Thus the two paths are disjoint
ordinary simple paths, giving one canonical pattern with
$a_1=2k$.

Every vertex other than an original path endpoint has positive
indegree and outdegree. Hence there are at most four degree-one
vertices, at most two vertices with zero indegree, and at most
two with zero outdegree. Since $\sum_x\deg x=2(a+r)$ and
$a+r-v=q-2$, the degree sums give
\[
 \begin{aligned}
 \sum_{\deg x\ge3}(\deg x-2)
 &=2q-4+\#\{x:\deg x=1\}\le2q,\\
 \sum_{\deg^\pm(x)\ge2}(\deg^\pm(x)-1)
 &=q-2+\#\{x:\deg^\pm(x)=0\}\le q.
 \end{aligned}
\]
The first identity gives at most $2q$ vertices of degree at least three, with total
degree at most $6q$. The second bounds the total indegree at
vertices of indegree at least two by $2q$, and likewise for
outdegree. Including the at most four additional endpoints of degree at
most two, there are at most $2q+4$ branch vertices, with total
degree at most $6q+8$.

Every inconsistent edge uses one of these incoming
or outgoing half edges, so there are at most $4q$ such edges.
Also, every raw edge belongs to a remainder block; hence there are at most $rk\leq qk$ raw single edges.

\medskip
\noindent\emph{Counting and class weights.}
Use the encoding of Lemma~\ref{lem:mixed-count}.
The total degree of branch vertices is at most $6q+8$, so there are at most
$3q+4$ maximal non-branching paths. Reserve $3q+4$ triple slots for these paths,
as in Step~1 in the proof of Lemma~\ref{lem:mixed-count}. There are $3(3q+4)$ integers.

The two block pairs are recorded as in Step~2, and the two
blocks give $2+r\le4$ indexed tangle-free paths whose lengths
are determined by the block types and cuts. This requires $4$ integers.

Choose a spanning forest as in Lemma~\ref{lem:mixed-count}.
Since
\[
 \sum_H g(H)=q+c-2\le q,
\]
there are at most $q$ edges outside the forest. Each compressed
tangle-free path traverses every edge at most twice, so its
list of edges traversed outside the forest has at most $2q$
entries. It therefore uses at most $2q+3$ records: its initial
and final vertices, and the list length and entries, with $8q+12$ integers in total.

For half-edge identifications, we record one integer per incoming or outgoing support edge at each branch. The vertices of degree at least three have total degree at most $6q$. There are at most four additional original endpoints, each of degree at most two. Therefore,
\[
\sum_{x\ {\rm branch}}\deg x\le 6q+8.
\]
Putting together,
these records use at most \[3(3q+4)+4+(8q+12)+(6q+8)\leq 36(q+1)\] integers.
The two paths have $2k$ steps in total, so the augmented
support has at most $2k$ edges and $2k+2$ vertices, and we may encode with
integers $\{0,\ldots,2k+2\}$. The encoding therefore gives
\[
 \#\{\gamma:q(\gamma)=q\}\le(2k+3)^{36(q+1)}.
\]

For the weights, apply the local degree calculation of
Lemma~\ref{lem:mixed-weight}. An ordinary vertex is defined to have degree
two and is not a branch or incident to a formal
edge. It has $r_x=h_x^+=h_x^-=1$ and contributes
$\beta_{t_x}$, where $r_x$ is its augmented indegree.
The number of ordinary vertices with $t_x=1$ is at most
$a_1+(k+4)q$, by the bounds on single edge counts above.
All other vertices are exceptional: they are branches or
degree-two vertices incident to formal edges. Their total degree is at
most $6q+8+4r\le 6q+16$.
The ordinary vertices have $r_x=1$, so
\[\sum_{x\ {\rm exceptional}}(1-r_x)=v-(a+r)=2-q.\]
Thus 
\[
 n^{-2}\prod_{x\ {\rm exceptional}}n^{1-r_x}
 =n^{-2}n^{2-q}=n^{-q}.
\]
Since $\sum_xt_x=2k$ and $2k=o(\sqrt N)$, the matching
estimate from Lemma~\ref{lem:matching} and the local degree bounds in Lemma~\ref{lem:beta}  give
\[
 |W(\gamma)|
 \le C2^{6q}\Delta^{6q+16}n^{-q}
       \rho^{2k-a_1-(k+4)q}
       \left(\frac{6k}{\sqrt{2n}}\right)^{a_1}.
\]
Multiplying by the count, and using
$\rho^{-1}\le\sqrt\Delta$ and $\Delta^k\le n^{1/8}$ as in
Lemma~\ref{lem:mixed-errors}, gives
\[
\begin{aligned}
 \rho^{-2k}
 \sum_{\substack{\gamma:q(\gamma)=q\\a_1(\gamma)=a_1}}|W(\gamma)|
 &\le Cn^{-\frac{15q}{16}-\frac{a_1}{2}}
       (4\Delta(k+1))^{44q+52+a_1}\le Cn^{-\frac{15q}{16}-\frac{a_1}{2}}
       (4\Delta(k+1))^{52(q+a_1+1)}.
\end{aligned}
\]

\medskip
\noindent\emph{Summation and the scalar power.}
Since $\Delta = n^{o(1)}$ and $k=O(\log n)$,
$(16\Delta(k+1))^{52}=n^{o(1)}$. Hence the two geometric
ratios
\[
 \frac{(16\Delta(k+1))^{52}}{n^{15/16}},
 \qquad
 \frac{(16\Delta(k+1))^{52}}{\sqrt n}
\]
are at most $1/2$ for sufficiently large $n$.
The case $(q,a_1)=(0,0)$ is absent, so summing over $q+a_1\ge1$
by the same calculation as in Lemma~\ref{lem:mixed-errors} yields
\[
 \frac{\E z_k^2}{\rho^{2k}}
 \le\frac{8C(4\Delta(\ell+1))^{104}}{\sqrt n}
 =n^{-1/2+o(1)}.
\]
By Markov's inequality,
\[
 \Pp\{|z_k|> n^{-1/8} \rho^k\}
 \le\frac{n^{1/4} \E z_k^2}{\rho^{2k}}=o(1).
\]
This proves \eqref{eq:scalar}.
\end{proof}

\subsection{Proof of Theorem~\ref{lem:CLT}}

Lemma~\ref{lem:local}, applied with $\alpha=1/8$, gives
$\Pp(\cH_\ell^c)=o(1)$.
On $\cH_\ell$, Lemma~\ref{lem:identification} gives
\[
 \frac{X_{\ell,t}-z_{\ell,t}}{\sqrt{t}\,\rho^{t}}
 =\frac{\tr P^{t}-1}{\sqrt{t}\,\rho^{t}}.
\]
By Lemma \ref{lem:scalar}, $z_{\ell,t} / \rho^{t}\to 0$ in probability. Then Lemma~\ref{lem:moments} and Slutsky's
theorem give a Gaussian limit for the left-hand side.
The equality on $\cH_\ell$ 
then transfers this limit to the trace:
\begin{align}
     \frac{\tr P^{t}-1}{\sqrt{t}\,\rho^{t}}
 \Longrightarrow\mathcal N(0,1).
\end{align}

\subsection*{AI statements} The authors developed the proof ideas, performed the calculations, and wrote the manuscript. Large language models were used solely for polishing the language and creating the figures.

\subsection*{Acknowledgments} The authors thank Jorge Garza-Vargas for helpful discussion during the early stage of this project. C.B. has started this work during his stay at the Institute for Advanced Study. He acknowledges
support from the James D. Wolfensohn Fund and from French ANR grant Tagada (ANR-25-CE40-
5672).  Y.Z. was partially supported by the Simons Grant MPS-TSM-00013944 and NSF DMS-2606337. Part of the work was done while C.B. and Y.Z. were visiting the Simons Institute for the Theory of Computing during the Spectral Theory Beyond Graphs program in Fall 2026. 

\bibliographystyle{plain}
\bibliography{ref}
\end{document}